\documentclass[reqno]{amsart}

\usepackage[all]{xy}
\usepackage{graphicx}

\usepackage{amssymb}
\usepackage{amsmath}
\usepackage{mathrsfs}
\usepackage{epsfig}
\usepackage{amscd}
\usepackage{graphicx, color}

\def\E{\ifmmode{\mathbb E}\else{$\mathbb E$}\fi} %natural numbers
\def\N{\ifmmode{\mathbb N}\else{$\mathbb N$}\fi} %natural numbers%
\def\R{\ifmmode{\mathbb R}\else{$\mathbb R$}\fi} %real numbers
\def\Q{\ifmmode{\mathbb Q}\else{$\mathbb Q$}\fi} %rational numbers
\def\C{\ifmmode{\mathbb C}\else{$\mathbb C$}\fi} %complex numbers
\def\H{\ifmmode{\mathbb H}\else{$\mathbb H$}\fi} %complex numbers
\def\Z{\ifmmode{\mathbb Z}\else{$\mathbb Z$}\fi} %integers
\def\P{\ifmmode{\mathbb P}\else{$\mathbb P$}\fi} %real numbers
\def\T{\ifmmode{\mathbb T}\else{$\mathbb T$}\fi} %real numbers
\def\SS{\ifmmode{\mathbb S}\else{$\mathbb S$}\fi} %real numbers
\def\DD{\ifmmode{\mathbb D}\else{$\mathbb D$}\fi} %real numbers

\newcommand{\del}{\partial}

\newcommand{\ben}{\begin{enumerate}}
\newcommand{\een}{\end{enumerate}}
\newcommand{\be}{\begin{equation}}
\newcommand{\ee}{\end{equation}}
\newcommand{\bea}{\begin{eqnarray}}
\newcommand{\eea}{\end{eqnarray}}
\newcommand{\beastar}{\begin{eqnarray*}}
\newcommand{\eeastar}{\end{eqnarray*}}
\newcommand{\bc}{\begin{center}}
\newcommand{\ec}{\end{center}}

\theoremstyle{theorem}
\newtheorem{thm}{Theorem}[section]
\newtheorem{cor}[thm]{Corollary}
\newtheorem{lem}[thm]{Lemma}
\newtheorem{prop}[thm]{Proposition}

\theoremstyle{definition}
\newtheorem{defn}[thm]{Definition}
\newtheorem{rem}[thm]{Remark}
\newtheorem{ques}[thm]{Question}

\newtheorem*{thm*}{Theorem}

\numberwithin{equation}{section}

\def\R{{\mathbb R}}

\def\E{{\mathbb E}}
\def\Z{{\mathbb Z}}
\def\C{{\mathbb C}}
\def\R{{\mathbb R}}
\def\P{{\mathbb P}}

\def\N{{\mathbb N}}

\def\11{{\mathbb I}}

\def\delbar{{\overline \partial}}
\def\dudtau{{\frac{\del u}{\del \tau}}}
\def\dudt{{\frac{\del u}{\del t}}}

\def\C{\mathbb{C}}
\def\Z{\mathbb{Z}}

\def\T{\mathbb{T}}

\def\Q{\mathbb{Q}}

\def\E{\ifmmode{\mathbb E}\else{$\mathbb E$}\fi} %natural numbers
\def\N{\ifmmode{\mathbb N}\else{$\mathbb N$}\fi} %natural numbers
\def\R{\ifmmode{\mathbb R}\else{$\mathbb R$}\fi} %real numbers
\def\Q{\ifmmode{\mathbb Q}\else{$\mathbb Q$}\fi} %rational numbers
\def\C{\ifmmode{\mathbb C}\else{$\mathbb C$}\fi} %complex numbers
\def\H{\ifmmode{\mathbb H}\else{$\mathbb H$}\fi} %complex numbers
\def\Z{\ifmmode{\mathbb Z}\else{$\mathbb Z$}\fi} %integers
\def\P{\ifmmode{\mathbb P}\else{$\mathbb P$}\fi} %real numbers
\def\SS{\ifmmode{\mathbb S}\else{$\mathbb S$}\fi} %real numbers
\def\DD{\ifmmode{\mathbb D}\else{$\mathbb D$}\fi} %real numbers

\def\R{{\mathbb R}}

\def\E{{\mathbb E}}
\def\Z{{\mathbb Z}}
\def\C{{\mathbb C}}
\def\R{{\mathbb R}}

\def\N{{\mathbb N}}
\def\LL{{\mathcal L}}

\def\delbar{{\overline \partial}}

 \def\ep{\epsilon}

\def\CA{{\mathcal A}}

\def\CJ{{\mathcal J}}

\def\CL{{\mathcal L}}

\def\darr#1{\raise1.5ex\hbox{$\leftrightarrow$}
\mkern-16.5mu #1}

\def\roughly#1{\raise.3ex\hbox{$#1$\kern-.75em
\lower1ex\hbox{$\sim$}}}

\def\opname#1{\mathop{\kern0pt{\rm #1}}\nolimits}

\def\span{\operatorname{span}}

\begin{document}
\quad \vskip1.375truein

\def\mq{\mathfrak{q}}
\def\mp{\mathfrak{p}}
\def\mH{\mathfrak{H}}
\def\mh{\mathfrak{h}}
\def\ma{\mathfrak{a}}
\def\ms{\mathfrak{s}}
\def\mm{\mathfrak{m}}
\def\mn{\mathfrak{n}}
\def\mz{\mathfrak{z}}
\def\mw{\mathfrak{w}}
\def\Hoch{{\tt Hoch}}
\def\mt{\mathfrak{t}}
\def\ml{\mathfrak{l}}
\def\mT{\mathfrak{T}}
\def\mL{\mathfrak{L}}
\def\mg{\mathfrak{g}}
\def\md{\mathfrak{d}}
\def\mr{\mathfrak{r}}

\title[Lagrangian Floer theory: energy and $C^0$-estimates]{
Unwrapped continuation invariance in Lagrangian Floer theory: energy and $C^0$ estimates}

\author{Yong-Geun Oh}
\address{Department of Mathematics, University of Wisconsin-Madison, Madison, WI, 53706
\& National Institute for Mathematical Sciences,
Daeduk Boulevard 628, Yuseong-gu, Daejeon,
305-340 Korea}
\email{oh@math.wisc.edu}
\thanks{The author is partially supported by the NSF grant DMS-0904197}

\date{October 4, 2009}

\begin{abstract}
We consider pairs of Lagrangian submanifolds $(L_0,L), \,
(L_1, L)$ belonging to the class of Lagrangian submanifolds with
\emph{conic} ends on \emph{Weinstein manifolds}.
The main purpose of the present paper is to define
a canonical chain map $h_\CL: CF(L_0,L) \to CF(L_1,L)$ of Lagrangian
Floer complex inducing an isomorphism in homology,
under the Hamiltonian isotopy $\CL=\{L_s\}_{0 \leq
s\leq 1}$ generated by \emph{conic} Hamiltonian functions such that the
intersections $L \cap L_s$ do not escape to infinity. The main
ingredients of the proof is an a priori bound
for general isotopy of the energy \emph{quadratic} at infinity and
a $C^0$-bound for the \emph{$C^1$-small} isotopy $\CL = \{L_s\}$, for the
associated pseudo-holomorphic map equations with \emph{moving}
Lagrangian boundary induced by a conic Hamiltonian isotopy.
For the Lagrangian submanifolds with \emph{asymptotically conic} ends,
we construct a natural homomorphism $h_\LL: HF(L_0,L) \to HF(L_1,L)$
for which the corresponding chain map may \emph{not} necessarily exist.
%
%This provides a more conventional construction of the chain isomorphism which
%replaces the sophisticated method using the Lagrangian cobordism via
%the machinery of \cite{kasturi-oh1,kasturi-oh2} whose details were
%only outlined in \cite{oh:gokova}.
\end{abstract}


\maketitle

\tableofcontents

\section{Introduction}
\label{sec:intro}

Floer \cite{floer:morse} invented the Floer homology $HF(L_0, L_1)$ of the pair $(L_0,
L_1)$ of \emph{compact} Lagrangian submanifolds with suitable topological restrictions
on \emph{compact} symplectic manifolds $(M,\omega )$
(or more generally on $M$ with bounded geometry). He defined this by
considering the (generalized) Cauchy-Riemann equation
\be\label{eq:CR}
\begin{cases}
\dudtau +J\dudt =0 \\
u(\tau ,0)\in L_0,\;\; u(\tau ,1)\in L_1
\end{cases}
\ee
for a map $u:\R\times [0,1]\to M$ and a one-parameter family of
almost complex structures $J = \{J_t\}_{0 \leq t \leq 1}$.

One crucial property of $HF(L_0,L_1)$ for applications to the problems in
symplectic topology is the invariance property under the Hamiltonian
deformations of the pair. Floer's original proof [F1] considers the case
where $L_1=\phi^1_H(L_0)$ and $\pi_1 (P,L_0)=\{ e\}$ where $\phi^1_H:M\to
P$ is the time-one map of the Hamiltonian flow of the function $H:[0,t] \times M \to\R$,
and involves some combinatorial study of the changes
occurring to the boundary operators when a (generic) degenerate
intersection occurs between the pairs during the deformations.
But this construction is \emph{not} canonical.

Motivated by the approach taken in \cite{floer:cmp} for Hamiltonian
diffeomorphisms and by the consideration of `moving' Lagrangian boundary
in \cite{hofer1}, the present author \cite{oh:cpam} used a variant of
\eqref{eq:CR}
\be\label{eq:CRchain}
\begin{cases}
\dudtau + J \dudt=0 \\
u(\tau ,0)\in L, \; u(\tau ,1)\in L_{\chi (\tau )}
\end{cases}
\ee
to construct a canonical chain map
where $\chi :\R\to [0,1]$ is a monotonically increasing function with
$\chi (-\infty )=0$ and $\chi (+\infty )=1$, i.e., a pseudo-holomorphic
equation with `moving' Lagrangian boundary condition.
\par
This approach works as long as the pair $(L_0, L_1)$ is compact on
geometrically bounded $M$. One of the crucial ingredients to work
with the moduli spaces of solutions of \eqref{eq:CR} and of
\eqref{eq:CRchain} on non-compact $M$ is the a priori energy bound and the
$C^0$ bound of the solutions $u$. It turns out that for the pair
$(L_0, L_1)$ which are \emph{conic} on \emph{Weinstein} manifold
$M$, i.e., of the type $[1,\infty) \times R\subset \R \times Q$ at
the end with $R$ a Legendrian submanifold in a contact manifold $Q$ a
simple application of the strong maximum principle (and the maximum
principle) enables one to prove this $C^0$-bound for \eqref{eq:CR}.
(See \cite{EHS} for the first such application, and \cite{oh:jdg},
\cite{seidel:biased} later for such applications.)
\emph{When one has an energy bound}, the monotonicity type argument
\cite{sikorav} can also be used to prove the $C^0$ bound.
(See \cite{nadler-zaslow,nadler}.)
\par
However the author later observed \cite{kasturi-oh2,oh:gokova} that application of
(strong) maximum principle for the equation \eqref{eq:CRchain} fails to hold in general
but works only for the isotopy of Lagrangian submanifold $(L_{0,s}, L_{1,s})$
\emph{in some monotone direction}. In \cite{oh:jdg}, a spectral invariant
$\rho(H;S)$ was assigned to a submanifold $S \subset N$ by considering the Floer
homology $HF(\nu^*S, o_N;H)$ and the author just stated that the invariant
is continuous under the $C^1$-topology of submanifolds. However to prove this
continuity, the author should have studied how the filtration of the Floer
chain complex behaves under the chain map
$$
h_{\alpha\beta}:HF(\nu^*S_\alpha,o_N;H) \to HF(\nu^*S_\beta,o_N;H)
$$
and especially should have examined existence of the chain map under the
Hamiltonian isotopy. However he presumed such a chain map exists based on his experience
with the compact case in \cite{oh:jdg} but later noticed that this
presumption is ill-founded because the analysis to construct the chain map for \emph{noncompact} Lagrangian
submanifolds such as $\nu^*S$ requires the energy and the $C^0$-bounds which
do not easily follow unlike the compact case. Since he could not
prove this $C^0$-bound in the framework of \cite{oh:cpam}, the author
outlined a different scheme in \cite{oh:gokova} of constructing a
natural chain map by considering a suspension cobordism of the isotopy $\{L_s\}$ on the suspension
$T^*\R \times M$ of $M$ and applying the geometric constructions developed in
\cite{kasturi-oh1,kasturi-oh2}. Although the author has little doubt
that this approach can be completed as explained therein, the
details have not been carried out yet.
Recently there have been many literature studying the Floer theory of
noncompact Lagrangian submanifolds in relation to the mirror symmetry. (See
\cite{HIV}, \cite{nadler-zaslow}, \cite{nadler}, \cite{seidel:book}, \cite{FSS1,FSS2} and etc.)
This continuation invariance is an important ingredient especially in relation
to the application to problems of symplectic topology. However its proof
in the literature is somewhat murky in details, at least to the present author.

The main purpose of the present paper is to rectify the status of this matter by
providing a more conventional construction of the chain map under the isotopy
$\CL = \{L_s\}$ of exact Lagrangian submanifolds with \emph{conic ends}
whose definition we will make precise later in  section \ref{sec:weinstein}. We will do this by establishing the required
a priori energy bound of the solutions of \eqref{eq:CRchain}
for general \emph{conic} Hamiltonian isotopies of such Lagrangian submanifolds
in general Weinstein manifold, but the $C^0$ bound only for a \emph{$C^1$-small} isotopy,
which will be enough to construct the \emph{adiabatic} chain map as in \cite{milinko-oh}.
A similar statement was proven in the cotangent bundle in \cite{abb-schwarz}
for the case of \emph{superlinear} growth:
but the energy bound in this general context of Weinstein manifolds outside
the cotangent bundle seems to hold only for the Hamiltonian isotopy of \emph{conic} Hamiltonian
but not for the Hamiltonians with superlinear growth, \emph{unless} the Hamiltonian
satisfies various restrictions on its growth or on the sign behavior.
Our energy bound enables us to apply the bubbling and compactness argument
to obtain the required $C^0$-estimates for a $C^1$-small isotopy.
(See \cite{abb-schwarz} for a study of the energy estimates and $C^0$-estimates
for the chain map on the cotangent bundle. Their proof of $C^0$-estimates
does not seem to generalize to the current context.)

The class of smooth conic Lagrangian submanifolds
includes all conormals of submanifolds, more generally, (smoothed) conormals of \emph{standard pairs}
treated in \cite{kasturi-oh1,kasturi-oh2} or the (smoothed) micro-supports of constructible
functions \cite{KS}. The class of asymptotically conic ones includes the \emph{standard (or costandard) Lagrangians}
of the cotangent bundle \cite{nadler-zaslow,nadler} and the \emph{Lagrangian thimbles}
\cite{HIV,seidel:book,FSS1} and their deformations generated by conic Hamiltonians.

The main objective of this paper is to give a complete proof of the following theorem.
This theorem was stated in \cite{oh:gokova} but its proof was only outlined by
a different method therein.

\begin{thm}\label{thm:conic} Let $M$ be a Weinstein manifold and assume $(L_0, L)$, $(L_1, L)$ are
transverse pairs of either compact or conic exact
Lagrangian submanifolds or their mixture.
Suppose that $\CL = \{L_s\}$ is a Hamiltonian isotopy generated by
a Hamiltonian function of conic type and that the set
$
\bigcup_{s \in [0,1]} L_s \cap L
$
is compact. Then there exists a natural chain map
$$
h_{\CL}: CF(L_0,L) \to CF(L_1,L)
$$
which induces an isomorphism in homology.
\end{thm}

The main theorem boils down to proving an a priori energy bound and
a $C^0$-bound for the solutions \eqref{eq:CRchain}.

After we take care of the conic case, we consider \emph{asymptotically conic
ones}. See section \ref{sec:weinstein} for
the precise definition of Lagrangian submanifolds with asymptotically conic ends.
Once we have made this definition precise, we can provide a canonical
procedure of approximating an asymptotically conic Lagrangian submanifolds by
conic ones, we can define the natural isomorphism
$$
h_{\LL}: HF(L_0,L) \to HF(L_1,L); \quad \LL = \{L_s\}_{s \in [0,1]}
$$
as the inverse limit of the homomorphisms
$$
h_{\LL_R}: HF(L_{0,R},L_R) \to HF(L_{1,R},L_R)
$$
of the approximations $(L_{0,R},L_R)$, $(L_{(1,R)},L_R)$ as $R \to \infty$. We denote this limit
homomorphism by $h_{\LL}$.

\begin{thm}\label{thm:asympconic} Let $M$ be a Weinstein manifold and assume $(L_0,L)$, $(L_0,L_1)$
be as in Theorem \ref{thm:conic} except that the hypothesis of being conic replaced
by being asymptotically conic. Suppose $\LL = \{L_s\}$ is a Hamiltonian isotopy
generated by a conic Hamiltonian function and the intersection $\cup_{s \in [0,1]} L_s \cap L$
is compact. Then there exists a natural isomorphism
$$
h_\LL: HF(L_0,L) \to HF(L_1,L).
$$
\end{thm}

We would like to point out that we do not know whether
there exists a relevant chain map
$CF(L_0,L) \to CF(L_1,L)$ that induces this isomorphism. This somewhat resembles the situation
in \cite{kasturi-oh1,kasturi-oh2} in relation to constructing a
Floer Fary functor associated to a (compact) exact Lagrangian submanifold
in the cotangent bundle.

\begin{ques} Does there exist a chain map
$$
h_\LL: CF(L_0,L) \to CF(L_1,L)
$$
that induces the above mentioned homomorphism $HF(L_0,L) \to HF(L_1,L)$?
\end{ques}

It is the author's impression that resolving this question is an
important one because many literature do not seem to examine this issue carefully when
they treat the Floer homology of asymptotically conic Lagrangian submanifolds
such as standard Lagrangians \cite{nadler-zaslow,nadler} or
Lagrangian thimbles \cite{HIV,seidel:book,FSS1}. We begin to suspect that the answer
to this question is negative in general. If this is indeed so, one may have to pay some
caution in the chain level Floer theory of asymptotically conic
Lagrangian submanifolds used in the literature.

Organization of the contents of the paper is now in order.
After a brief review of some standard background materials concerning
Weinstein manifolds and construction of chain maps in Floer theory for
the compact Lagrangian submanifolds from \cite{oh:cpam},
section \ref{sec:chainmap}-\ref{sec:wrapup} treat the case of
pairs $(L,L_0)$ of conic Lagrangian submanifolds.

Section \ref{sec:energybound} proves an energy bound which
\emph{depends on} the conic threshold $R_\LL$ of the isotopy $\LL$.
One novelty of our proof of the energy bound (for the chain map)
is our usage of the energy induced by the symplectic form that has the `quadratic growth'
at infinity, not the `linear growth' arising from the symplectization.
While this bound certainly implies the energy bound
arising from the latter, we have not been able to directly
prove the latter bound using the metric associated to the
symplectic form on the symplectization. (See Remark \ref{rem:quadmetric}
for the reason why.)

Section \ref{sec:derivbound} proves a pointwise derivative bound.
Our proof of the pointwise derivative bound uses a proof by contradiction which
exploits the \emph{convexity at infinity} of the triple $(M,\omega, J)$ and
the special geometry of (asymptotically) conic Lagrangian submanifolds on Weinstein manifolds.
Here it is important to use the Courant-Lebesgue lemma together with
the monotonicity formula as used in \cite{oh:removal}.

Section \ref{sec:C0bound} considers the $C^0$-estimates for the case where
\emph{the Hamiltonian isotopy is sufficiently small}. Under this $C^1$-smallness
of the isotopy, we obtain a uniform $C^0$-bound
by using the energy bound and applying the maximum principle and the monotonicity
formula. \emph{While the $C^0$-bound itself does not depend on the conic threshold, the
above mentioned $C^1$-smallness depends on it, which is the reason why
we cannot construct a chain map for the asymptotically conic pairs}.

Section \ref{sec:wrapup} then concludes construction of the chain map for
the conic Lagrangian pairs by partitioning the given Hamiltonian isotopy into
small ones for which the above bounds, especially the $C^0$ bound, can be applied.

Finally section \ref{sec:asympconic} explains how we can take care of the
asymptotically conic pairs. Here we prove a Darboux-type theorem on the
asymptotic Lagrangian submanifolds in a neighborhood of a conic Lagrangian
submanifold, and construct an explicit family of approximations thereof by
the conic Lagrangian submanifolds. Then we take the inverse limit of the
continuation homomorphisms in $HF$ for the approximations and obtain
the continuation homomorphism for the given asymptotic conic pairs.

In the last section of this paper, we will explain how one can
extend the similar energy and $C^0$ bounds to the cases of general
\emph{non-exact} Lagrangian submanifolds of conic type on
non-compact symplectic manifolds with conic ends,
after taking into account of the bubbling phenomena. This kind of
study will be necessary for studying the invariance property of
invariants extracted from the `unwrapped' Fukaya type category
generated by such non-compact Lagrangian submanifolds. Here we use
the term `unwrapped' to contrast the type of end behavior considered
in this paper to those `wrapped' version considered in
\cite{abb-schwarz}, \cite{abou-seidel}, where a super linear growth
condition of the Hamiltonian functions is imposed. The same kind of
$C^0$-estimate is used in the author's proof of Weinstein's
conjecture on symplectically fillable contact manifolds \cite{oh:weinstein}.

We would like to thank Rezazadegan for asking a question on this
continuation invariance which has motivated us to write out a complete proof
thereof. We also thank Abbondandolo for several useful e-mail
communications in which he kindly explained details of the $C^0$-estimates
presented for the case of cotangent bundle in \cite{abb-schwarz}.

\section{Weinstein manifolds and Lagrangian submanifolds of conic type}
\label{sec:weinstein}

In this section, we recall from \cite{EG} that a symplectic manifold
$(M,\omega)$ is called \emph{convex at infinity} if it carries a
vector field $X$ which is {\it complete symplectically dilating} at
infinity: A vector field $X$ is complete symplectically dilating if
the flow $\{\phi^t\}$ of X is complete and satisfies $(\phi^t)^*\omega =
e^t \omega$. We assume that $(M,\omega)$ allows an exhausting
pluri-subharmonic function at infinity. Following \cite{EG}, we call
such a manifold {\it Weinstein} (at infinity).

We choose $\varphi$ an exhausting pluri-subharmonic function with respect to a tame almost
complex structure $J$. We also assume that $J$ is invariant under
the flow of $X$ outside a compact set. Then the level set
$\varphi^{-1}(S)$ for sufficiently large $S$ carries the induced
contact structure (in fact a $CR$-structure) on it, and the flow map
$$
\phi: (0,\infty) \times \varphi^{-1}(S) \to
M \cap \varphi^{-1}([S,\infty)) ; \, (s,y) \mapsto \phi^s(y)
$$
defines a diffeomorphism. In this coordinates, we have $\varphi = s$ and
$\phi^s$ is nothing but the translation maps $\phi^s(S,y) = (S+s,y)$.
With this said, the definition of a Weinstein manifold
can be given in terms of this contact manifold $Q = \varphi^{-1}(S)$
which is now in order.

We choose a contact form $\lambda$ of $Q$ so that
\begin{enumerate}
\item $\omega|_Q = d\lambda$ for some one-form $\lambda$
and $\xi = \ker \lambda$
\item the orientation of $Q$ defined by $\lambda \wedge
(d\lambda)^{n-1}$ coincides with the boundary orientation of $Q =
\partial W$.
\end{enumerate}
When this latter condition holds, the contact manifold $(Q,\xi)$ is called
\emph{strongly symplectically fillable}.

By the symplectic neighborhood theorem, one can choose
a function $r$ in a collar neighborhood $U_\delta$ of $\del W = Q$
in $W$ such that
\bea\label{eq:Udelta}
U_\delta & \cong & (1-\delta,1] \times Q \nonumber \\
\omega & = & d(r \pi^*\lambda) \quad \mbox{on } \, U_\delta \eea for
the projection $\pi : M^{\geq R} \to Q$. We then consider the
cylinder $(1-\delta, \infty) \times Q$ and form the union
\be\label{eq:decomp} \widehat W = W \# (1 -\delta, \infty) \times Q \ee along
the strip $(1-\delta, 1] \times Q \cong U_\delta$. The symplectic
form $\omega$ naturally extends to $\widehat W: = M$ by gluing it with $d(r
\pi^*\lambda)$ on $(1-\delta, \infty) \times Q$. We denote the corresponding symplectic
from by
$$
\widehat \omega = \omega \# d(r\pi^*\lambda).
$$
By definition, every
Weinstein manifold has this decomposition by identifying $\varphi = \log
r$, or equivalently $r = e^{\varphi}$ in terms of the identification
$M^{\geq R} \cong [0,\infty) \times \varphi^{-1}(\ln R)$. In terms of
the chosen contact form $\lambda$ on $Q$ and the projection $\pi :
M^{\geq R} \to Q$, the symplectic form $\omega$ has the expression
$$
\omega = d(r \pi^*\lambda) = d(e^{s} \wedge
\pi^*\lambda), \quad s = \varphi.
$$
On $Q$, the Reeb vector field $X_\lambda$ associated to the contact
form $\lambda$ is the unique vector field satisfying
\be\label{eq:Liouville} X \rfloor \lambda = 1, \quad X \rfloor
d\lambda = 0.
\ee
Therefore the tangent bundle $TQ$ has the splitting $TQ =
\span\{X_\lambda\} \oplus \xi$. We denote by
$$
\pi_\lambda: TQ \to \xi
$$
the corresponding projection.
\par
We call $(s,y)$ the cylindrical coordinates and $(r,y)$ the cone
coordinates. In the case of $M = \C^n\setminus \{0\} \cong
(0,\infty) \times S^{2n-1}$, $(\sqrt{r},y)$ with $y \in S^{2n-1}$ is
nothing but the standard polar coordinates of $\C^n \setminus \{0\}$.
On the other hand if we write
$$
T^*N \setminus \{0\} \cong S^1(T^*N) \times \R_+
$$
then $r = |p|$ for the canonical coordinates $(q,p)$ of $T^*N
\setminus \{0\} \subset T^*N$.
On the cylinder $[0, \infty) \times Q \subset (-\infty,\infty) \times Q$,
we have the natural splitting
$$
TM \cong \R\cdot \frac{\del}{\del r} \oplus TQ \cong \span\left\{X_\lambda, \frac{\del}{\del s}\right\} \oplus \xi
\cong \R^2 \oplus \xi.
$$
We denote by $\widetilde X_\lambda$ the unique vector field on $M$ which is
invariant under the translation, tangent to the level sets of $r$
and projected to $X_\lambda$. When there is no danger of confusion, we will
sometimes just denote it by $X_\lambda$.

Now we describe a special family of almost complex structure adapted to
the given cylindrical structure of $M$.

\begin{defn} An almost complex structure $J$ on $(M,\omega)$ is called
\emph{$\lambda$-contact type} if it is split into
$$
J = j \oplus J_\xi : TM \cong \R^2 \oplus \xi \to TM \cong \R^2 \oplus \xi
$$
where $J|_\xi$ is compatible to $d\lambda|_{\xi}$ and $j: \R^2 \to \R^2$
maps $\frac{\del}{\del s}$ to $X_\lambda$.
\end{defn}

We denote by $g_Q$ the metric on $Q$ compatible to the contact form $\lambda$ and
the endomorphism $J_\xi$ defined by
$$
g_Q(h,k) = g_{(\lambda,J_\xi)}: = \lambda(h)\lambda(k) + d\lambda(h, J_\xi k).
$$
We will also need to consider the associated cylindrical metric on $(1,\infty) \times Q$
is given by
$$
g^{cl} = ds^2 + g_Q.
$$
For our purpose, we will need to consider a family of symplectic forms
to which the given $J$ is compatible and their associated metrics.

We consider a positive smooth function $\psi: (0, \infty) \to \R_+$ that satisfies
\be
\psi(r) = \begin{cases} r \quad & \mbox{for }\, 0 < r \leq 1 \\
r^2 & \mbox{for } \, r \geq 3
\end{cases}
\label{eq:psi}
\ee
and $\psi$ satisfies $r \leq \psi(r) \leq r^2$ and $\psi''(r),\, \psi'(r) > 0$
for $ 1 \leq r \leq 3$, and define the two from
$$
\omega_\psi: = d(\psi(r) \pi^*\lambda).
$$
For any $\lambda$-contact $J$, the $J$-compatible metric associated to $\omega_\psi$ is expressed as
\be\label{eq:gvarphi}
g_{(\psi,J)} =  \frac{\psi'(r)}{r}dr^2 + r \psi'(r) \lambda \otimes \lambda + \psi(r) g_Q(\pi_\lambda(\cdot), \pi_\lambda(\cdot)).
\ee
In particular we have the formula for the metric $g_{(\omega,J)}$ compatible to
$\omega$ 
\be\label{eq:gomega}
g_{(\omega,J)} = \frac{1}{r}(dr^2 + r^2 g_Q).
\ee
We denote  $g_{lin,J}$ and $g_{\operatorname{quad},J}$ to be $g_{(\psi,J)}$ with $\psi$
replaced by $r$ or $r^2$, respectively.  Then we have
$$
g_{(\omega,J)} \leq g_{(\psi,J)}.
$$
We will use both $g^{cl}$ and $g_{(\psi,J)}$ in a suitable way for our purpose.
The following obvious fact is one reason why we use the `quadratic' metric

\begin{lem}\label{lem:drnorm} Consider the metric $g_{(\psi,J)}$. Then we have
$$
|dr|_{g_{(\psi,J)}} = 1
$$
on $[3, \infty) \times Q$.
\end{lem}

\begin{defn} We call $H: M \to \R$ is
\emph{of conic type (at the end)}) if it is homogeneous of order 1, i.e., if
it satisfies
$$
m_\gamma^*dH = \gamma dH
$$
on $[R_0,\infty) \times Q$ for some $R_0$ where $m_\gamma$ is the
multiplication map $r \mapsto \gamma r$.
\end{defn}
In the cylindrical coordinates it satisfies
$$
\tau_{s_0}^*dH = e^{s_0}dH
$$
where $\tau_{s_0}$ is the translation $s \mapsto s + s_0$ .
We would like to compare this choice of conic Hamiltonians to the \emph{superlinear}
ones used in the \emph{wrapped} version of Floer theory studied in
\cite{abb-schwarz}, \cite{abou-seidel} and others:
\begin{defn}\label{defn:contacttype}
A function $H : M \to \R$ is called \emph{of contact type (at
the end)}, if it has the form
$$
H = h(e^S), \quad S(s,y) = s - f(y)
$$
at the end $M^{\geq R} = \varphi^{-1}([S, \infty))$ with $R = e^S$ for a positive
smooth function $h: \R_+ \to \R_+$ with $h' \geq  0$ and a function $f: Q \to \R$.
\end{defn}
In this context, they also put a condition corresponding to
$$
\lim_{r \to \infty}(rh'(r) - h(r)) = \infty
$$
which is responsible for the term `wrapped'. (See \cite{viterbo} for the first
usage of such conditions in relation to the Floer theory.)

Examples of Weinstein manifolds include cotangent bundles of compact
manifolds. A typical Hamiltonian of conic type is given as follows: consider an isotopy $\{ S^s\}_{0\le s\le
1}$ of submanifolds $S^s\subset  N$.  Consider the
corresponding deformation of the conormal bundles
$$
\{ \nu^*S^s\}_{0\le s\le 1}\subset T^*N.
 $$
This deformation is realized by the $s$-dependent Hamiltonian
\be\label{eq:conormalH} H(s, q,p)=\langle p, X_s(q)\rangle
\ee
where $X_s$ is the vector field realizing the isotopy $\{ S^s\}$ i.e.\
$X_s={d\over ds}\big| S^s$. Certainly, this family of Hamiltonians is \emph{neither} compactly supported
\emph{nor} of contact type in the sense of Definition \ref{defn:contacttype} nor
satisfies any sign condition. This is one of the author's motivation to
carefully examine the question of continuation invariance of Floer homology
$HF(\nu^*S, o_N)$ under the isotopy of submanifold $S \subset N$.

Next we introduce the notion of Lagrangian submanifolds with conic
(respectively, asymptotically conic) end.
We will just simply call such Lagrangian submanifolds conic
(respectively, asymptotically conic) Lagrangian
submanifolds.

\begin{defn} A Lagrangian submanifold $L$ is said to \emph{have conic end}, if
there exists some $R_0 > 0$ such that $L \cap M^{\geq R_0}$ is
invariant under the dilatation, or more precisely if it holds
$$
m_\gamma (L \cap M^{\geq R_0}) \subset L \cap M^{\geq R_0}
$$
for all $\gamma \geq 1$. We call the minimum value of such $R_0$ the \emph{conic threshold}
of $L$ and denote it by $R_L$.
\end{defn}

The following proposition will be important later for our derivation of
the energy bound, which explains why the `conic' condition for $L$ is
important. We emphasize that this proposition does not hold even for
`asymptotically conic' Lagrangian submanifolds, which is the reason why
we have not been able to treat the latter class of Lagrangian submanifolds
as much as the conic ones.

\begin{prop}\label{prop:Lagwrtomegapsi} Let $L$ be a Lagrangian submanifold with respect to
$\widehat \omega =d \alpha$, $\alpha = r\pi^*\lambda$
that has conic end of its conic threshold $R_L$. Consider a function $\psi$ satisfying \eqref{eq:psi}.
Suppose in addition that $\psi(r) \equiv r$ for $1 \leq r \leq R_L$.
Then $L$ is also Lagrangian with respect to the symplectic form $\omega_\psi$, and
if $i_L^*(\alpha) = df$ on $L$, then we also have $i_L^*(\alpha_\psi) = df$ with
$\alpha_\psi = \psi(r) \pi^*\lambda$.
\end{prop}
\begin{proof} Since $L$ is conic on $M^{\geq R_L}$, it follows
$$
L \cap [R_L, \infty) \times Q = N \times [R_L,\infty)
$$
for some Legendrian submanifold $N \subset (Q,\xi)$. Therefore we have
$$
i_L^*(r \pi^*\lambda) =
r\circ i_L \, (\pi \circ i_L)^*\lambda = r\circ i_L \, i_N^*\lambda \equiv 0.
$$
on $M^{\geq R_L}$. Similar computation proves $i_L^*(\psi(r) \pi^*\lambda) = 0$
thereon.
The second statement follows since $df = 0$ on $L\cap M^{\geq R_L}$
and $\psi(r) \pi^*\lambda = r \pi^*\lambda$ on $L\cap M^{\leq R_L}$.
\end{proof}

It takes some preparation to provide a rigorous definition of
asymptotically conic end.
Intuitively, we say $L$ \emph{has asymptotically conic end} if
\be\label{eq:asympconic}
\lim_{\gamma \to \infty} m_\gamma^{-1}(L \cap M^{> \gamma R_0}) = (R_0, \infty) \times N^\infty
\ee
for some $R_0$. We call $N^\infty \subset Q$ the \emph{asymptotic boundary} of $L$
and denote it by $\del_\infty L = L^\infty$.
To give the rigorous definition, we need to make the topology of convergence
in \eqref{eq:asympconic} more precise. Denote $L \cap r^{-1}(R) = \{R \} \times
L(R)$ in the identification of $M^R = \{R\} \times Q$ i.e.,
$$
L(R) = \pi(L \cap r^{-1}(R)) \subset Q
$$
and consider the $C^k$-distance $\operatorname{dist}_{C^k}(L(R), N)$ for a
given compact $(n-1)$-dimensional submanifold $N \subset Q$.

\begin{defn} \label{defn:asympconic}
For a given $R_0$, we denote
\be\label{eq:asympconic}
\lim_{\gamma \to \infty} m_\gamma^{-1}(L \cap M^{\geq \gamma R_0}) = [R_0, \infty) \times N
\ee
if $\operatorname{dist}_{C^0}(L(R), N) \to 0$ uniformly over $R \in [R_0,\infty)$.
We call $N \subset Q$ the \emph{asymptotic boundary} of $L$
and denote it by $\del_\infty L = N$. We say $L$ has \emph{$C^k$-asymptotically conic end}
if $\operatorname{dist}_{C^k}(L(R), N) \to 0$ uniformly over $R \in [R_0,\infty)$.
\end{defn}

The following lemma is easy to prove.

\begin{lem} If a Lagrangian submanifold $L$ has $C^1$-asymptotically conic end,
then the asymptotic boundary $\del_\infty L$ is Legendrian.
\end{lem}

This gives rise to the following natural question of $C^0$-symplectic
geometry.

\begin{ques} Suppose that a Lagrangian submanifold $L$ has
$C^0$-asymptotically conic end and that the asymptotic boundary $\del_\infty L$
is $C^1$-submanifold of $Q$. Will $\del_\infty L$ be Legendrian?
\end{ques}

We would like to alert readers that the intersection $L \cap r^{-1}(R_0)$ is isotropic
but \emph{not} Legendrian on $\{R_0\} \times Q$ \emph{unless $L$ is conic}.

The following $C^0$-estimate for the Floer boundary operator
can be easily proven by a version of strong maximum
principle (See [EHS] for the details).

\begin{lem}\label{lem:maximum} Let $(M,\omega)$ be as above and
consider Lagrangian submanifolds $L,\, L_0$.
Let $j=\{J_t\}_{0\leq t \leq 1}$ be a family of almost complex structures
such that $J_t = J$ outside a compact set.
Suppose that $L_0 \cap L_1$ are compact and $L_i$'s are
transverse to the level set $r^{-1}(R_0)$ for $R_0 \geq 1$ and assume that
the intersections of them with $r^{-1}(R_0)$ are Legendrian in the
level set $r^{-1}(R_0)$ with respect to the induced contact structure.
Then the maximum $r\circ u$ cannot be achieved at any point
$(\tau,t)$ with $r(u(\tau,t)) = R_0$ for any solutions $u$
of \eqref{eq:CR}.
\end{lem}

But the strong maximum principle (nor the monotonicity formula)
cannot be applied to \eqref{eq:CRchain}
for the continuity equation in general.
It applies in some monotone direction of the isotopy $\CL = \{L_s\}$.
(See \cite{kasturi-oh2}, \cite{seidel:biased} for such remarks.)

\begin{rem}
We would like to take this chance to point out that
there is an error in the statement of Theorem 2.1 \cite{oh:gokova}
in that the hypothesis
\begin{itemize}
\item $L_i$'s are transverse to the level sets of $\psi$ at infinity
\end{itemize}
should be replaced by
\begin{itemize}
\item $L_i$'s are conic at infinity.
\end{itemize}
Strong maximum principle applies only under the latter assumption, but not
under the former. While one can apply the monotonicity argument instead to prove the
energy or the $C^0$ bound \emph{for the boundary map} as illustrated by
\cite{sikorav,nadler-zaslow} (in the cotangent bundle case), this does \emph{not} seem to
apply to \emph{the case of the chain map} which involves a \emph{moving} boundary.
In relation to this, it appears to the author that
the proof of continuation invariance presented in \cite{nadler-zaslow,nadler},
especially the one presented in the appendix of \cite{nadler},
is not satisfactory.
\end{rem}

\section{Action functional, Hamiltonian isotopies and Floer chain maps}
\label{sec:chainmap}

From now on, when we say $L$ is Lagrangian, it will always mean to be with respect to the
symplectic form $\widehat \omega = \omega \# d(r \pi^* \lambda)$ for a
given contact form $\lambda$ on $Q$.

In this section, we review the construction of the canonical chain map under
Hamiltonian isotopies introduced in \cite{oh:cpam} for the Lagrangian Floer homology
on \emph{closed} or more generally geometrically bounded, not necessarily
compact, symplectic manifolds. In this paper, we will
restrict our attention to the case of \emph{exact} Lagrangian submanifolds.

We first recall the corresponding action functional on the path space
of the pair two exact Lagrangian submanifolds $((L,f_L), (L_0,f_{L_0}))$ in an exact
symplectic manifold $(M, \omega)$ with $\omega = d\alpha$. Here
the functions $f_L:L \to \R$ and $f_{L_0}: L \to \R$ are chosen so that
$df_L = i_L^*\alpha$ and $df_{L_0} = i_{L_0}^*\alpha$. They are defined up to
addition of constants.
Then we consider the functional $\CA : \Omega(L_0,L) \to \R$ by
\be\label{eq:AA}
\CA(\gamma) = - \int \gamma^*\alpha - f_{L}(\gamma(1)) + f_{L_0}(\gamma(0)).
\ee

Let $J = \{J_t\}_{0 \leq t\leq 1}$ be a one-parameter family of compatible
almost complex structures. Equip $\Omega(L_0,L)$ with the $L^2$-metric
$\langle\langle \cdot, \cdot \rangle \rangle_J$ defined by
$$
\langle\langle \xi_1, \xi_2 \rangle \rangle_J = \int_0^1 g_{J_t}(\xi_1(t),\xi_2(t))\, dt
$$
where $g_{J_t}$ is the Riemannian metric on $M$ defined by $g_{J_t} = \omega(\cdot, J_t \cdot)$.
A simple computation shows
\be\label{eq:gradAA}
\operatorname{grad} \CA(\gamma) = J_t \dot \gamma
\ee
with respect to the compatible metric $g_J = \omega(\cdot, J\cdot)$ and so the negative $L^2$-gradient equation for
$\CA$ is given by
\be\label{eq:CRbdy}
\begin{cases}
\dudtau + J_t \dudt = 0\\
u(\tau,0) \in L_0, \, u(\tau,1) \in L.
\end{cases}
\ee
The following energy estimate is the crucial element for the analysis of the
moduli space of solutions of this equation.

\begin{prop} For any solution of \eqref{eq:CRbdy} with finite energy and bounded image,
$u(\tau) \to p^\pm \in L_0 \cap L$ as $\tau \to \pm \infty$
we have
\beastar
E_{\widehat \omega,J}(u) & = & \CA(\widehat p^-)- \CA(\widehat p^+) = f_{L_0}(p^+) - f_L(p^-)\\
& \leq & \max_{p,\, p' \in L_0 \cap L}(f_{L_0}(p) - f_L(p'))=: C_0
\eeastar
In particular if $L_0 \cap L$ is compact, we have uniform upper bound
for the energy independent of such solution $u$.
\end{prop}

We refer to \cite{sikorav} for the details of the proof of the following monotonicity
formula.

\begin{prop}[Monotonicity]\label{prop:monotonicity} Consider the $(M,L,L_0, g_{(\psi,J)})$.
Then there exists $r_0 > 0$ depending only on the bounds of the curvature of $g_{(\psi,J)}$,
the injectivity radius, and of the norms of second fundamental forms of $L, \, L_0$ such that
the following holds:
Let $u: \R \times [0,1] \cong D^2 \setminus\{0,1\} \to M$ be a
$J$-holomorphic map
$$
u(\pm \infty) \in L_0 \cap L, \, u(\tau,0) \in L_0, \, u(\tau,1) \in L_1.
$$
Assume $u(S) \subset B(x,r)$ for some domain $S \subset \R \times [0,1]$ with
$u(\del S\setminus \R \times \{0,1\}) \subset \del B(x,r)$ and $x \in u(S)$. Then there
exists $C_5 > 0$ with $r_0$ as above such that
\be\label{eq:monotonicity}
\operatorname{Area}_{g_{(\psi,J)}}(u(S)) \geq C_5 r^2
\ee
for any $0 < r \leq r_0$.
\end{prop}

Using the above energy bound and the monotonicity formula, we immediately
obtain the following $C^0$ bound too. (See \cite{nadler-zaslow} for such a proof.)
For the readers' convenience, we include its proof.

\begin{cor} Suppose $J_t \equiv J_0$, $J_0$ is of $\lambda$-contact type.
Suppose that $L, \, L_0$ are $C^2$-asymptotically conical Lagrangian submanifolds.
Then for for any solution with finite energy and bounded image we have
$$
\operatorname{Im} u \subset M^{\leq R_0}
$$
where $R_0$ is at least smaller than
$$
R_{L_0 \cap L} + \frac{2C_0}{C_5 r_0}
$$
where $R_{L_0 \cap L}$ is the infimum of $R$ such that $L_0 \cap L \subset M^{\leq R}$.
\end{cor}
\begin{proof} The proof is based on the  monotonicity formula.
Due to the assumption on the $C^2$-asymptotic conical condition, the norms of
the second fundamental forms of $L, \, L_0$ are bounded and so we can apply
the monotonicity formula.

We choose $R=R(u) > 0$ depending on $u$ such that
$$
\operatorname{Im } u \subset M^{\leq R}
$$
for all finite energy solution $u$ of \eqref{eq:dudtauXFdudt} with bounded image.
Now we want to estimate the infimum of such $R$'s.

We fix the constant $r_0> 0$ appearing in Proposition \ref{prop:monotonicity}
above. Then we can find at least $\frac{R-R_{L_0\cap L}}{2r_0}$
disjoint balls $B(x,r_0)$ with $x \in \operatorname{Im}u$. Therefore the above
monotonicity formula implies
$$
C_5 r_0^2 \cdot \frac{R-R_{L_0 \cap L}}{2r_0} \leq C_0
$$
Therefore we obtain
$$
R \leq R_{L_0 \cap L} + \frac{2C_0}{C_5 r_0}
$$
which finishes the proof.
\end{proof}

This energy estimate is the starting point of the construction of Floer
boundary map in general.

\begin{rem} When the Lagrangian submanifolds are conic, not just asymptotically conic,
one can apply also the strong maximum principle to prove a stronger $C^0$ bound.
(See \cite{EHS}, \cite{oh:gokova}.)
But the strong maximum principle does not apply to general
\emph{asymptotically conic} Lagrangian submanifolds, though.
\end{rem}

Next we consider the isotopy $\CL = \{L_s\}$ connecting $L_0$ to $L_1$,
and the pseudo-holomorphic equation with \emph{moving} boundary condition
\be\label{eq:moving}
\begin{cases}
{\partial u\over\partial\tau} + J_{t} {\partial
u\over\partial
t}=0 \\
u(\tau ,0)\in L_{\chi(\tau)},\;\; u(\tau ,1)\in L
\end{cases}
\ee
for a cut-off function $\chi:\R \to [0,1]$ given by
\be\label{eq:chi}
\chi(\tau) = \begin{cases} 1 \quad & \mbox{for } \, \tau \geq 1\\
0 \quad & \mbox{for } \, \tau \leq 0.
\end{cases}
\ee
We will fix this cut-off function once and for all and not change it
and so this dependence on $\chi$ will be sometimes suppressed in our exposition.

According to the scheme given in \cite{oh:cpam}, the chain map
$$
h_\CL:HF(L_0,L)\rightarrow HF(L_1,L)
$$
is defined by considering the isolated solutions of \eqref{eq:moving}. This
is called the \emph{geometric} version of the Floer equation.
Although the construction of the chain map using the \emph{moving} Lagrangian boundary
is carried out on compact closed $(M,\omega)$ in \cite{oh:cpam}, it also applies
to the case of  non-compact $(M,\omega )$
\emph{as long as} one can prove the energy bound and the $C^0$-bound
for the solutions of the equation. Unlike the equation \eqref{eq:CRbdy}
proving either bound turns out to be non-trivial due to the moving
boundary condition. The proofs of these two bounds are indeed precisely
the main technical results of the present paper.

Another way of defining the chain map, which is sometimes easier to handle when
obtaining various geometric estimates, is to consider the \emph{dynamical} version which
uses the fixed boundary condition but with the equation perturbed by
the Hamiltonian vector field. The explanation of this transformation is now in order.

Suppose the isotopy $\{L_s\}$ is achieved by the Hamiltonian isotopy
$\phi=\{\phi^s\}$ generated by a $s$-dependent Hamiltonian functions $F_1 = F_1(s,x)$ of
conic type for $0 \leq s \leq 1$. In particular we have $L_1 = \phi_{F_1}^1(L_0)$.
We interpolate $F_1$ and the 0 function via a two-parameter family of conic
Hamiltonians $F:[0,1]\times [0,1] \times M \to \R$ such that
$$
F(s,1,x) = F_1(s,x), \quad F(s,0,x) = 0
$$
e.g., we can choose $F(s,t,x) = \rho(t) F_1(s,x)$ for a surjective increasing
function $\rho:[0,1] \to [0,1]$ such that $\rho(0) = 0 , \rho(1)= 1$ and
$\rho$ is constant near $t =0, \, 1$.

Then we denote $F_t(s,x): = F(s,t,x)$ and $\phi(s,t)$ the Hamiltonian flow of
$F_t$ in the direction of $s$, i.e.,
$$
\phi(s,t) = \phi_{F_t}^s.
$$
We note that the isotopy $\phi(s,1)$ is the original isotopy $\phi^s$
realizing the isotopy $\{L_s\}$ from $L_0$ to $L_1$ while $\phi(0,1) \equiv id$ and
so it fixes $L$ at $t = 0$.

We consider the composition map
$$
\widetilde u(s,t) = (\phi(s,t))^{-1}\circ u(s,t).
$$
One can easily check that $\widetilde u$ satisfies the perturbed
Cauchy-Riemann equation \be\label{eq:XFXHchi}
\begin{cases}
\left({\partial \widetilde u\over\partial\tau} -
X_{F^\chi}(\widetilde u) \right)
 + (\phi(\chi(\tau),t)^*J \left({\partial \widetilde u\over\partial
t} - X_{H_t}(\widetilde u)\right) =0  \\
\widetilde u(\tau ,0)\in L_0,\;\; \widetilde u(\tau ,1)\in L
\end{cases}
\ee
where $X_{F^s}$ (respectively $X_{H_t}$) is the Hamiltonian vector field generated by
$F_t$ in the direction of $s$ (respectively, $F^s$ in the
direction of $t$).

Although we will not emphasize it in this paper
it is in general necessary to consider  families of
diffeomorphisms, Hamiltonians and almost complex structures depending
on the domain parameter of the pseudo-holomorphic curves. See \cite{abou-seidel,seidel:biased}
and other papers of Seidel for an extensive systematic treatment of
such a usage. The method of our proofs of energy and $C^0$-bound
also applies to this general situation which we will not mention since
we will not use it in this paper. However this will be needed to
prove the continuation invariance in the context of Fukaya category in
general. See section \ref{sec:nonexact} for some discussion on this point.

\section{Energy bound}
\label{sec:energybound}

Let $\{L_s\}_{0\le s\le 1}$ be a Hamiltonian isotopy associated to a
Hamiltonian
$$
F_1=F_1(s,x): [0,1] \times M \rightarrow \R
$$
of conic type at the end. If we let $R_{\LL} = \max\{R_{L_s}\}$,
we can choose $F_1$ so that it becomes conic on $R^{\geq R_\LL}$.
By reparameterizing the isotopy $L_s$ near $s = 0, \, 1$ so that it stays constant
near there, we may assume that
\be\label{eq:0atend}
F_1(s,\cdot) \equiv 0 \quad \mbox{near } \, s = 0, \, 1.
\ee
\emph{We would like to note that this adjustment will play some role in our proof
of a priori energy bound later.}

We then define their elongations
$F_1^{\chi}: \R \times [0,1] \times M \to \R$ by reparameterizing
$s$ by $s = \chi(\tau)$ and setting $F_1^\chi(\tau,t,x) = F_1(\chi(\tau),x)$. Then we have
\be\label{eq:0atinfty}
F_1^{\chi(\tau)} \equiv 0 \quad \mbox{ for }\, \tau \leq 0\,  \mbox{ or }\, \tau \geq 1.
\ee

Motivated by the discussion in the previous section
we extend the Hamiltonian isotopy $\phi_{F_1}^s$ to a 2-parameter family
$$
\phi: [0,1] \times [0,1] \to Ham(M,\omega)
$$
Hamiltonian diffeomorphisms such that
$$
\phi(s,0) \equiv id, \quad \phi(s,1) = \phi_{F_1}^s.
$$
In particular this isotopy at $t = 1$ moves $L_0$ to $L_1$ in the $s$-direction
and it does not move $L$ in the $s$-direction at $t=0$.

By a suitable $t$-reparameterization, we may assume that $\phi(s,\cdot)$ is
constant in $t$ near $t = 0,\, 1 \in [0,1]$. We denote the Hamiltonian generating the vector fields
by
$$
X_{F^s}(s,t,x) = \frac{\del \phi}{\del s} \circ \phi^{-1},
\, \quad X_{F_t}(s,t,x) =  \frac{\del \phi}{\del t} \circ \phi^{-1}
$$
respectively. Here we denote
$$
F^s(t,x) = F(s,t,x) = F_t(s,x)
$$
depending on the direction of $s$ or $t$ in which we take the Hamiltonian
vector fields. Then the composition
$$
\widetilde u(s,t) = \phi(s,t)^{-1}\circ u(s,t)
$$
satisfies \eqref{eq:XFXHchi}.

Similarly we consider a two-parameter family $J =\{J_{(s,t)}\}$ of $\lambda$-contact almost complex
structures
$
J:[0,1] \times [0,1]\rightarrow \CJ_\omega
$
such that
\be\label{eq:Jondelsquare}
J_{(s,t)} \equiv J_0 \quad \mbox{for }\, (s,t) \in \del [0,1]^2
\ee
and its elongation $J^\chi = \{J_{\chi(\tau), t}\}: \R \times [0,1] \to \CJ_\omega$.
\par
Consider two pairs of Lagrangian  submanifolds
$(L_0,L)$ and $(L_1,L)$ for which $L_0,\, L_1$ and $L$ are either
compact or conic at the end. Since other cases are easier to handle,
we will assume that all the pairs are non-compact and conic at the end.
We will always assume that
$\bigcup_{s \in [0,1]} L \cap L_s$ is a compact subset of $M$.
This condition means that during the isotopy $\{L_s\}$ from $s =0$ and $s =1$, there
is no intersection escaping to infinity.

We consider the following function $\psi$ given by
\be\label{eq:psiRR0}
\psi(r) = \begin{cases} r \quad & \mbox{for } \, 1 \leq r \leq R \\
r^2 + (R-R^2) \quad & \mbox{for } \, r \geq R +1
\end{cases}
\ee
for $R = R_\LL$ and suitably interpolated in between so that $1 \leq \psi'(r) \leq 2(R_\LL +1)$.
We note that
$$
\psi(r) \geq r, \quad \frac{\psi'(r)}{r} \geq \frac{1}{r}, \quad \psi(r) \leq r^2
$$
and so we still have
$$
g_{(\widehat \omega,J)} \leq g_{(\omega_\psi,J)} \leq g_{\operatorname{quad},J}.
$$

Motivated by the discussion in the last section, we now
consider a smooth solution $u:\R\times [0,1]\rightarrow T^*M$ of
the dynamical version of \eqref{eq:CRchain},
\be\label{eq:FHCR}
\begin{cases}
\left(\dudtau - X_{F^{\chi(\tau)}}(u)\right) + J^{\chi}\left(\dudt - X_{F_t}(u)\right) =0 \\
u(\tau ,0) \in L, \quad u(\tau,1)\in L_0.
\end{cases}
\ee
We emphasize that this equation has a \emph{fixed} boundary condition.

Now we consider the symplectic forms $\omega_{\psi} = d(\psi(r) \pi^*\lambda)$ and denote the associated
family of metrics
$$
g_{(\psi,J^\chi)} = \{ g_{(\psi,J_{\chi(\tau),t})} \}_{(\tau,t) \in \R \times [0,1]}.
$$
By Proposition \ref{prop:Lagwrtomegapsi}, $L, \, L_0$ are still Lagrangian with respect to
$\omega_\psi$. We define the energy of a smooth
map $u: \R \times [0,1] \to M$ by
$$
E_{(g_{(\psi,J^\chi)},F^\chi)}(u) = \frac{1}{2}
\int\int\left(\left|\dudtau - X_{F^{\chi(\tau)}}(u)\right|_{g_{(\psi,J^\chi)}}^2
+\left|\dudt - X_{F_t}(u)\right|_{g_{(\psi,J^\chi)}}^2\right) \, dt\, d\tau
$$
which is reduced to
\beastar
E_{(g_{(\psi,J^\chi)},F^\chi)}(u) & = & \int\int\left|\dudtau - X_{F^{\chi(\tau)}}(u)\right|_{g_{(\psi,J^\chi)}}^2\, dt\, d\tau\\
& = & \int\int\left|\dudt - X_{F_t}(u)\right|_{g_{(\psi,J^\chi)}}^2\, dt\, d\tau
\eeastar
for a solution $u$ of \eqref{eq:FHCR}.
\par
To make the most useful statement of the theorem, we introduce a few
geometric invariants of the pair $(L_0,L)$.

First due to the assumption that $(L_0,L)$ is transverse and $L_0 \cap L$ is compact,
there are a finite number of points in it. We consider the action functional
$\CA$ associated to the pair $((L_0,f_{L_0}),(L,f_L))$ and the
differences
$\CA(\widehat x_\alpha) - \CA(\widehat x_\beta)$ for all possible
$x_\alpha, \, x_\beta \in L_0 \cap L$
and define the \emph{action width}
\be\label{eq:actionwidth}
w_\CA(L_0,L) = \max\{|\CA(\widehat x_\alpha) - \CA(\widehat x_\beta)|
\mid x_\alpha, \, x_\beta \in L_0 \cap L \, \}.
\ee
(In fact this can be written in term of the functions $f_{L_\alpha}$ and $f_{L_\beta}$
by the formula \eqref{eq:AA}. Since this explicit form is not used in this paper,
we do not write down the precise expression for this.)

Next consider a conic function $F$ such that $F$ has the form $F(s,t,(r,y)) =
f(s,t,y) r$ on $M^{\geq R_\LL}$ or equivalently
$$
F(s,t,x) = r(x) f (s,t,\pi(x))
$$
where $f$ is a function defined on $[0,1]^2 \times Q \to \R$.
Here we set $R_L = \infty$ if $L$ is compact.
We define the quantities
\be\label{eq:e-int}
\ep_{(F,J)}^{in}(R_\LL) = \max_{(s,t,x) \in [0,1]^2 \times M^{\leq R_\LL}} |dF(s,t,x)|_{g_{(\psi,J_t)}}
\ee
and
\bea\label{eq:e-out}
\ep_{(F,J)}^{out;0}& = & \max_{(s,t,y) \in [0,1]^2 \times Q} \{f(s,t,y) \, \mid \, (s,t) \in [0,1]^2, \, y \in Q\,\}
\nonumber \\
\ep_{(F,J)}^{out;1} & = & \max_{(s,t,y) \in [0,1]^2 \times Q} \{|df_{s,t}|_{g}(y)\, \mid
(s,t) \in [0,1]^2, \, y \in Q\, \}\nonumber \\
\ep_{(F,J)}^{out} & = & \ep_{(F,J)}^{out;0} + \ep_{(F,J)}^{out;1}.
\eea
The latter does not depend on $R_\LL$ but depends only on the
asymptotic limit of the function $F$ of conic type.
We then denote
\be\label{eq:eF}
\ep_{(F,J)} = \max\{\ep_{(F,J)}^{in}(R_\LL), \ep_{(F,J)}^{out}\}
\ee
From now on, we always assume $R_\LL \geq 3$ in this section.

The following a priori bound is the main theorem in this section.
We would like to emphasize that we cannot apply the maximum principle to obtain an
a priori $C^0$-estimates due to the Hamiltonian perturbation terms $X_{F^\chi}$ in the equation.
Hence we do not have the energy bound either. We also remark that
we do \emph{not} impose any sort of monotonicity on $F$ used such as in
\cite{floer-hofer}, \cite{kasturi-oh1,kasturi-oh2}, \cite{abou-seidel} which
would facilitate study of the energy and the $C^0$-estimates.

We note that if we assume the intersection set
$$
\operatorname{Int} \LL: = \bigcup_{s \in [0,1]}  L \cap L_s
$$
is compact, it is easy to see that it must be contained in $M^{\leq R_\LL}$:
Otherwise the conic nature of the triple $(M;L_0,L)$ and the isotopy $\LL$
would make the intersection set itself become conic and cannot be compact.

We will prove an energy estimate by comparing
the energy $E_{g_{(\psi,J^\chi)},F}$ with the change of action functional
$\CA_\psi$ in terms of the symplectic form $\omega_\psi$.
Since $L, \, \, L_0$ are Lagrangian with respect to $\widehat\omega$, which are a priori
not necessarily Lagrangian with respect to $\omega_\psi$, we cannot do this comparison
unless we ensure they also become Lagrangian with respect to $\omega_\psi$. This is
how Proposition \ref{prop:Lagwrtomegapsi} enters in our energy estimate and the difference
in our treatment of the conic and asymptotically conic cases arises.

\begin{thm}\label{thm:energybounds} Assume that $J$ and
$F$ as above. Then there exists  $C_1 = C_1(F,J)> 0$ for
$i=1, \, 2$ such that every finite energy solution $u$ of
\eqref{eq:FHCR} with
$$
u(-\infty) = \widehat x^-, \, u(\infty)=\widehat x^+
$$
for some $x^\pm \in L_0 \cap L$ satisfies
\be\label{eq:energy}
E_{(g_\psi^\chi,F)}(u) \leq C_1
\ee
where $C_1$ is at least smaller than
\be\label{eq:boundC1}
2\left(3R_\LL \ep_{(F,J)}^2 + w_\CA(L_0,L)\right).
\ee
\end{thm}

Since  $g_{(\omega_\psi,J)} \geq g_{(\widehat\omega,J)}$, we obviously have

\begin{cor}\label{cor:boundforomegahat}
We have the bound
$$
E_{(g_{(\widehat\omega,J)},F)}(u) \leq 2(3R_\LL\, \ep_{(F,J)}^2 + w_\CA(L_0,L)).
$$
\end{cor}

We find that the proof of this energy bound is somewhat curious in that the proof of the
energy estimate is carried out with the metric associated to the symplectic form
$$
\omega_\psi = d(\psi(r) \pi^*\lambda)
$$
not with the one induced by the symplectization form $\widehat\omega$. We
have not been able to prove the latter estimate directly using the metric $g_{(\widehat \omega,J)}$.

\begin{proof}[Proof of Theorem \ref{thm:energybounds}]
We calculate the variation of $\CA_\psi(u(\tau))$ where
$\CA_\psi: \Omega(L_0,L) \to \R$ is the action functional \eqref{eq:AA} with respect
to the symplectic form $\omega_\psi$, i.e.,
$$
\CA_\psi(\gamma) = - \int \gamma^*(\psi (r) \pi^*\lambda) - f_L(\gamma(1)) + f_{L_0}(\gamma(0)).
$$
Proposition \ref{prop:Lagwrtomegapsi} implies that $L, \, L_0$ are still Lagrangian
with respect to $\omega_\psi$ and so we have the gradient of $\CA_\psi$
with respect to the $L^2$-metric $\langle\langle \cdot, \cdot \rangle\rangle_J$ on
$\Omega(L_0,L)$ given by
\be\label{eq:gradCA}
\operatorname{grad}_J \CA_\psi(z) = J \dot z.
\ee
Therefore the negative $L^2$-gradient flow of this functional
is again the standard
$$
\begin{cases}
\dudtau + J \dudt = 0\\
u(\tau,0) \in L_0, \, u(\tau,1) \in L.
\end{cases}
$$
We would like to emphasize that this functional $\CA_\psi$ is the \emph{same} one
for all $\tau$.

We compute
\beastar
\frac{d}{d\tau}\CA_\psi(u(\tau)) & = & d\CA_\psi(u(\tau))\left(\dudtau\right)
= \int_0^1 \left\langle \dudtau, J^\chi\dudt\right \rangle_{g_{(\psi,J^\chi)}} \, dt \\
& = & \int_0^1 \left\langle \dudtau, J^\chi\left(\dudt - X_{F_t}(u)\right) \right \rangle_{g_{(\psi,J^\chi)}} \, dt
+ \int_0^1 \left\langle \dudtau, J^\chi  X_{F_t}(u) \right \rangle_{g_{(\psi,J^\chi)}} \, dt\\
& = & \int_0^1 \left\langle X_{F^{\chi(\tau)}}(u)-J^\chi\left(\dudt- X_{F_t}(u)\right), J^\chi\left(\dudt - X_{F_t}(u)\right) \right \rangle_{g_{(\psi,J^\chi)}} \, dt\\
&{}& \quad + \int_0^1 \left\langle X_{F^{\chi(\tau)}}(u)-J^\chi \left(\dudt- X_{F_t}(u)\right), J^\chi X_{F_t}(u) \right\rangle_{g_{(\psi,J^\chi)}} \, dt\\
& = & - \int_0^1 \left|\dudt - X_{F_t}(u)\right|_{g_{(\psi,J^\chi)}}^2 \, dt +
\int_0^1 \left\langle X_{F^{\chi(\tau)}}(u), J^\chi  X_{F_t}(u)\right\rangle_{g_{(\psi,J^\chi)}}\, dt \\
&{}& \quad + \int_0^1 \left\langle X_{F^{\chi(\tau)}}(u), J^\chi \left(\dudt - X_{F_t}(u)\right) \right\rangle_{g_{(\psi,J^\chi)}}\, dt \\
&{}& \quad - \int_0^1 \left\langle J^\chi\left(\dudt - X_{F_t}(u)\right), J^\chi X_{F_t}(u) \right\rangle_{g_{(\psi,J^\chi)}} \, dt.
\eeastar
In the above derivation, we use the following facts in order:
\begin{enumerate}
\item We use the fact the function $\CA_\psi$ does not explicitly depend on $\tau$ for
the first equality,
\item we use  \eqref{eq:gradCA} for the second equality,
\item and we use the equation \eqref{eq:XFXHchi} for the fourth equality.
\end{enumerate}
Therefore we obtain
\bea\label{eq:|dutau|}
\int_0^1
&{}& \left|\dudt - X_{F_t}(u)\right|^2_{J^{\chi(\tau)}_t} \, dt
= - \frac{d}{d\tau}\CA_\psi(u(\tau))+
\int_0^1 \left\langle X_{F^{\chi(\tau)}}(u), J^\chi  X_{F_t}(u)\right\rangle_{g_{(\psi,J^\chi)}}\, dt \nonumber \\
&{}& \qquad+\int_0^1 \left\langle X_{F^{\chi(\tau)}}(u), J^\chi \left(\dudt - X_{F_t}(u)\right) \right\rangle_{g_{(\psi,J^\chi)}}\, dt\\
&{}& \qquad- \int_0^1 \left\langle J^\chi \left(\dudt - X_{F_t}(u)\right), J^\chi X_{F_t}(u) \right \rangle_{g_{(\psi,J^\chi)}} \, dt.
\label{eq:JchiXFt}
\eea
Integrating this over $\tau \in \R$, we obtain
\bea\label{eq:|dudt-XF|}
&{}& \int_{-\infty}^\infty \int_0^1
\left|\dudt - X_{F_t}(u)\right|^2_{J^{\chi(\tau)}_t} \, dt\, d\tau \\
&= & \CA_\psi(\widehat x^\alpha)- \CA_\psi(\widehat x^\beta) +
\int_{-\infty}^\infty \int_0^1 \left\langle X_{F^{\chi(\tau)}}(u), J^\chi  X_{F_t}(u)\right\rangle_{g_{(\psi,J^\chi)}}\, dt\, d\tau\nonumber\\
& {} & \quad +
\int_{-\infty}^\infty \int_0^1 \left\langle X_{F^{\chi(\tau)}}(u), J^\chi \left(\dudt - X_{F_t}(u)\right) \right\rangle_{g_{(\psi,J^\chi)}}\,
dt\, d\tau
\label{eq:intXFchi} \\
&{}& \quad - \int_{-\infty}^\infty \int_0^1 \left\langle J^\chi \left(\dudt - X_{F_t}(u)\right), J^\chi X_{F_t}(u) \right \rangle_{g_{(\psi,J^\chi)}}
\, dt\, d\tau
\label{eq:intXFt}\\
\nonumber
\eea
under the assumption
$$
E_{(g_\psi^\chi,F)}(u) = \int_{-\infty}^\infty \int_0^1
\left|\dudt - X_{F_t}(u)\right|^2_{J^{\chi(\tau)}_t} \, dt\,d\tau < \infty.
$$
It remains to estimate the terms \eqref{eq:intXFchi} and \eqref{eq:intXFt}.
\emph{Here is the place
where the requirement for $F$ to be conic and the usage of
the metric $g_{(\psi,J^\chi)}$, especially its quadratic nature at infinity,
play essential roles in our estimate.}

First we estimate \eqref{eq:intXFchi}.
Recall that
$$
F^{\chi(\tau)}  =  0 \, \mbox{ if $\tau\le 0$ or
if $\tau\ge 1$}
$$
from \eqref{eq:0atend}. Therefore have
\beastar
&{}& \int_{-\infty}^\infty \int_0^1 \left\langle X_{F^{\chi(\tau)}}(u), J^\chi\left(\dudt -X_{F_t}(u)\right)\right \rangle_{g_{(\psi,J^\chi)}} \, dt\, d\tau\\
& = & \int_0^1\int_0^1 \left\langle X_{F^{\chi(\tau)}}(u), J^\chi\left(\dudt -X_{F_t}(u)\right)\right \rangle_{g_{(\psi,J^\chi)}} \, dt\, d\tau.
\eeastar
Hence
\bea\label{eq:dFchiu0101}
&{}& \left|\int_{-\infty}^\infty \int_0^1 \left\langle X_{F^{\chi(\tau)}}(u), J^\chi\left(\dudt -X_{F_t}(u)\right)\right \rangle_{g_{(\psi,J^\chi)}} \, dt\, d\tau\right|\nonumber\\
& \leq & \int_0^1\int_0^1 \left|\left\langle X_{F^{\chi(\tau)}}(u), J^\chi\left(\dudt -X_{F_t}(u)\right)\right \rangle_{g_{(\psi,J^\chi)}}\right| \, dt\, d\tau\
\eea
When $r(u(\tau,t)) \leq R_\LL$, we have
\bea\label{eq:on<R0}
& {}& \left| \left\langle X_{F^{\chi(\tau)}}(u), J^\chi\left(\dudt -X_{F_t}(u)\right)\right \rangle_{g_{(\psi,J^\chi)}}\right|
\nonumber \\
& \leq & |X_{F^{\chi(\tau)}}(u)|_{g_{(\psi,J^\chi)}} \left|\dudt-X_{F_t}(u) \right|_{g_{(\psi,J^\chi)}}\nonumber\\
& = & |dF^\chi(u)|_{g_{(\psi,J^\chi)}}
\left|\dudt-X_{F_t}(u) \right|_{g_{(\psi,J^\chi)}} 
\leq \ep^{in}_{(F,J)}(R_\LL)
\left|\dudt-X_{F_t}(u) \right|_{g_{(\psi,J^\chi)}}
\eea
It follows by definition of the compatible metric $g_{(\psi,J^\chi)}$ and its
dual that
\be\label{eq:|XF|=|dF|}
|X_{F^s}(x)|_{g_{(\psi,J^\chi)}} = |dF^s(x)|_{g_{(\psi,J^\chi)}}
\ee
where the norms are taken in $TM$ and in $T^*M$ with respect to the metric
$g_{(\psi,J^\chi)}$ on $TM$  and its dual metric on $T^*M$ respectively.

On the other hand when $r(u(\tau,t)) \geq R_\LL$, the homogeneity of $F$ implies
$$
F(s,t,(r,y)) = f(s,t,y) r
$$
in coordinates $(s,t,(r,y)) \in [0,1]^2 \times [R_\LL,\infty) \times Q$ for some
function $f: [0,1]^t \times Q \to \R$.

\begin{lem} Let $r(x) \geq 1$. Then we have
$$
|dF^\chi(\tau,t,x)|_{g_{(\psi,J^\chi)}} \leq \sqrt{R_\LL}\cdot \ep_{(F,J)}^{out}.
$$
\end{lem}
\begin{proof}
We have
$$
dF^\chi = f^\chi\circ \pi dr + r d(f^\chi\circ \pi).
$$
We obtain
\bea\label{eq:|f||df|}
|f^\chi\circ \pi dr|_{g_{(\psi,J^\chi)}}
& \leq & |f^\chi|_\infty |dr|_{g_{(\psi,J^\chi)}}\leq \sqrt{R_\LL}\cdot |f^\chi|_\infty  \nonumber\\
|rd(f^\chi \circ \pi)|_{g_{(\psi,J^\chi)}} & = &|r df^\chi d \pi|_{g_{(\psi,J^\chi)}}
\leq \sqrt{R_\LL} \cdot |df^\chi|_\infty
\eea
where we use the fact that $r \geq 3$ and so $\psi(r) \equiv r^2$
and the derivative of the map $\pi: \{r\} \times Q \to \{1\} \times Q \cong Q$
has the norms
$$
|d\pi(r,y)|_{g_{(\psi,J^\chi)}} = 1/\sqrt{r}, \quad
|dr|_{g_{(\psi,J^\chi)}} \equiv \sqrt{r}
$$
on the linear region $1 \leq r \leq R_\LL$ and
$$
|d\pi(r,y)|_{g_{(\psi,J^\chi)}} = 1, \quad
|dr|_{g_{(\psi,J^\chi)}} \equiv 1
$$
on the quadratic region $r \geq R_\LL+1$. See \eqref{eq:gomega}. This finishes the proof.
\end{proof}

Combining these with \eqref{eq:on<R0}, we obtain
\beastar
& {}& \left| \left\langle X_{F^{\chi(\tau)}}(u), J^\chi\left(\dudt -X_{F_t}(u)\right)\right \rangle_{g_{(\psi,J^\chi)}}\right| \\
& \leq & \max\left\{\ep_{(F,J)}^{in}(R_0),\sqrt{R_\LL}\left(\ep_{(F,J)}^{out,0}+ \ep_{(F,J)}^{out,1}\right)\right\}
\left|\dudt-X_{F_t}(u)\right|_{g_{(\psi,J^\chi)}}\\
& \leq & \sqrt{R_\LL}\cdot \ep_{(F,J)}\left|\dudt-X_{F_t}(u)\right|_{g_{(\psi,J^\chi)}}
\eeastar
since $R_\LL \geq 1$. Therefore we have
\bea\label{eq:eF1/4}
&{}& \int_0^1\int_0^1 \left|\left\langle X_{F^{\chi(\tau)}}(u), J^\chi\left(\dudt -X_{F_t}(u)\right)\right \rangle_{g_{(\psi,J^\chi)}}\right| \, dt\, d\tau
\nonumber\\
& \leq & \int_0^1\int_0^1 \sqrt{R_\LL}\, \ep_{(F,J)}\left|\dudt-X_{F_t}(u)\right|_{g_{(\psi,J^\chi)}}\, dt\, d\tau \nonumber \\
& \leq & \int_0^1\int_0^1 R_\LL \,\ep_{(F,J)}^2\, dt\, d\tau
+ \frac{1}{4}\int_0^1\int_0^1 \left|\dudt-X_{F_t}(u)\right|_{g_{(\psi,J^\chi)}}^2 \, dt \,d\tau
\nonumber\\
& \leq & R_\LL \, \ep_{(F,J)}^2+ \frac{1}{4} \int_{-\infty}^\infty \int_0^1 \left|\dudt-X_{F_t}(u)\right|_{g_{(\psi,J^\chi)}}^2 \, dt \,d\tau.
\eea
For the second inequality, we used the inequality $ 2ab \leq  a^2 +  b^2$.

We carry out the similar estimate for the term \eqref{eq:JchiXFt} and obtain
\beastar
&{}& \int_0^1\int_0^1 \left|\left\langle J^\chi\left(\dudt -X_{F_t}(u)\right), J^\chi X_{F_t}(u) \right \rangle_{g_{(\psi,J^\chi)}}\right| \, dt\, d\tau \\
&\leq & R_\LL \, \ep_{(F,J)}^2+ \frac{1}{4} \int_{-\infty}^\infty \int_0^1 \left|\dudt- X_{F_t}(u)\right|_{g_{(\psi,J^\chi)}}^2 \, dt \,d\tau.
\eeastar
\emph{This is the reason why we wanted to reduce the integral $\int_{-\infty}^\infty (\cdot) d\tau$ to
$\int_0^1(\cdot) \, d\tau$ where the adjustment made for $F$ so that \eqref{eq:0atend} holds
in the beginning of this section is essential.} Otherwise we would not have
been able to obtain this kind of inequality for the integral $\int_{-\infty}^\infty(\cdot) \, d\tau$.

On the other hand, we have
\beastar
&{}& \left|\int_{-\infty}^\infty \int_0^1 \left\langle X_{F^{\chi(\tau)}}(u), J^\chi X_{F_t}(u) \right\rangle_{g_{(\psi,J^\chi)}}\,
dt\, d\tau\right| \\
& = & \left|\int_0^1 \int_0^1 \left\langle X_{F^{\chi(\tau)}}(u), J^\chi X_{F_t}(u) \right\rangle_{g_{(\psi,J^\chi)}}\,
dt\, d\tau\right| \\
& \leq &\left(\sqrt{R_\LL}\, \ep_{(F,J)}\right)^2 = R_\LL\, \ep_{(F,J)}^2.
\eeastar
Substituting all these into \eqref{eq:|dudt-XF|}, we obtain
\beastar
&{}& \int_{-\infty}^\infty \int_0^1 \left|\dudt- X_{F_t}(u)\right|_{g_{(\psi,J^\chi)}}^2 \, dt \,d\tau\\
& \leq & \frac{1}{2}
\int_{-\infty}^\infty \int_0^1 \left|\dudt- X_{F_t}(u)\right|_{g_{(\psi,J^\chi)}}^2 \, dt \,d\tau
+  3R_\LL \, \ep_{(F,J)}^2 +
\CA_\psi(\widehat x_\alpha) - \CA_\psi(\widehat x_\beta)
\eeastar
and hence
$$
E_{(g_{(\psi,J^\chi)},F) }(u) \leq
2(3 R_\LL\, \ep_{(F,J)}^2+ w_\CA(L_0,L)).
$$
This finishes the proof with a constant $C_1$ which is smaller than at least
$$
2 (3 R_\LL \ep_{(F,J)}^2 + w_\CA(L_0,L)).
$$
\end{proof}

\begin{rem}\label{rem:quadmetric}
We like to note that in the above proof
the crucial inequalities we used to obtain the above energy bound
are \eqref{eq:|f||df|}. The choice of $\psi$ is driven by our attempt
to have these two bounds for the conic Hamiltonians:
The quadratic growth of the metric is the only such choice to obtain
the bound independent of $r$.
\end{rem}

\section{Pointwise derivative bound}
\label{sec:derivbound}

We first prove the following derivative bound.

\begin{thm}\label{prop:gradient} Let $g^{cl}$ be the metric on
$M$ defined by
$$
g^{cl} = \begin{cases} g_J\quad & \mbox{on } \, M^{\leq 1} = W\\
ds^2 + g_Q \quad & \mbox{on } \, M^{\geq 2}
\end{cases}
$$
interpolated between them on $[1, 2] \times Q$. There exists a constant $C_2 > 0$
depending only on $C_1$ and $(L,L_0;M, \ep_{(J,F)})$ but independent of $u$'s such that
\be\label{eq:maxderivbound}
\max_{z \in \R \times [0,1]} |du(z)|_{g^{cl}} \leq C_2
\ee
for all solution $u$'s of \eqref{eq:XFXHchi}.
\end{thm}
\begin{proof} We first note that we can arrange the interpolation made for the
definition of $g^{cl}$ above so that $g^{cl} \leq  g_{(\psi,J)}$. Therefore
\bea\label{dq:gcyleng}
&{}& \frac{1}{2} \int_{\R \times [0,1]} \left| \frac{\del u}{\del \tau}
- X_{F^{\chi(\tau)}}(u) \right|_{g^{cl}}^2 + \left|\frac{\del u}{\del t} - X_{F_t}(u)\right|_{g^{cl}}^2
\nonumber\\
&  = & E_{g^{cl}}(u)
\leq  E_{(g_\psi^\chi,F)}(u) \leq C_1.
\eea
Suppose to the contrary that there exists sequences of
$z_i = (\tau_i,t_i)$, and $u_i$'s such that $L_i: = |du_i(z_i)|_{g^{cl}} \nearrow \infty$.
By choosing a subsequence, we may assume $z_i \to z_\infty$
if the sequence $z_i =(\tau_i,t_i)$ is bounded, or $z_\infty \to \infty$
otherwise. In any case, the sequence $J_{z_i}$
will converge: it converges to $J_\infty$ for the first case
and  to $J_0$ for the second case by the definition $J_K= J^{\rho_K}$ in
the beginning of section \ref{sec:energybound}. We denote either of the limit almost
complex structures by $J_\infty$.

Applying Lemma 26 \cite{hofer2}, we can choose $\delta_i > 0$ and rechoose $z_i$
so that
\bea\label{eq:eizi}
&{}& \delta_i |du(z_i)|_{g^{cl}} \to \infty, \, \delta_i \to 0,\nonumber\\
&{}& |du(z)|_{g^{cl}} \leq 2 |du(z_i)|_{g^{cl}} \quad \mbox{if } |z -z_i| \leq \delta_i.
\eea
We will consider two cases separately.

First assume that $u_i(z_i) \in M^{\leq R}$ for all $i$ for some $R >0$.

In this case, we apply the standard rescaling argument
$$
\widetilde u_i(z) = u_i \left(z_i + \frac{z}{L_i}\right).
$$
Using the conformal invariance of the harmonic energy, we compute
\beastar
&{}& \frac{1}{2} \int_{D^2(\delta_i L_i)} \left| \frac{\del \widetilde u_i}{\del \tau} \right|_{g^{cl}}^2
+ \left|\frac{\del \widetilde u_i}{\del t} \right|_{g^{cl}}^2\, dt\, d\tau\\
& = &\frac{1}{2} \int_{D_{z_i}(\delta_i)} \left| \frac{\del u_i}{\del \tau}
\right|_{g^{cl}}^2 + \left|\frac{\del u_i}{\del t}\right|_{g^{cl}}^2 \, dt\, d\tau\\
&\leq &
 \int_{D_{z_i}(\delta_i)} \left| \frac{\del u_i}{\del \tau}
- X_{F^\chi}(u_i) \right|_{g^{cl}}^2
 + \left|\frac{\del u_i}{\del t}- X_{F_t}(u)\right|_{g^{cl}}^2 \\
&{}& \quad
+ \int_{D_{z_i}(\delta_i)} |X_{F^\chi}(u_i)|_{g^{cl}}^2\, dt\, d\tau
+ \int_{D_{z_i}(\delta_i)} |X_{F_t}(u_i)|_{g^{cl}}^2\, dt\, d\tau \\
&\leq &
 \int_{D_{z_i}(\delta_i)} \left| \frac{\del u_i}{\del \tau}
- X_{F^\chi}(u_i) \right|_{g_{(\psi,J^\chi)}}^2
 + \left|\frac{\del u_i}{\del t}- X_{F_t}(u)\right|_{g_{(\psi,J^\chi)}}^2 \\
&{}& + \int_{D_{z_i}(\delta_i)} |X_{F^\chi}(u_i)|_{g_{\psi}^\chi}^2\, dt\, d\tau
+ \int_{D_{z_i}(\delta_i)} |X_{F_t}(u_i)|_{g_{(\psi,J^\chi)}}^2\, dt\, d\tau \\
& \leq & 2C_1 + 2 R_\LL \ep^2_{(F,J)}\cdot \pi \delta_i^2
\eeastar
for all  $i$'s.
We also have the derivative bound
\be\label{eq:derivbound}
|d\widetilde u_i(z)|_{g^{cl}} \leq 2
\ee
for all $z \in D^2(\ep_iK_i)$. %Write $\widetilde u_i = (a_i, \Theta_i)$.
And we have
\beastar
\max \left|\frac{\del u_i}{\del \tau} + J\frac{\del u_i}{\del t}\right|_{g_{(\psi,J^\chi)}} & \leq &
\max \left|\frac{\del u_i}{\del \tau} + J\frac{\del u_i}{\del t}\right|_{g_{(\psi,J^\chi)}}\\
& = & \max |X_{F^\chi}(u_i)|_{g_{(\psi,J^\chi)}} + \max |X_{F_t}(u_i)|_{g_{(\psi,J^\chi)}} \\
& \leq & 2 \sqrt{R_\LL}\cdot \ep_{(F,J)}.
\eeastar
Therefore we have
\be\label{eq:delbarto0} \max
\left|\delbar_J \widetilde u_i\right|_{g^{cl}} \leq \frac{1}{L_i}
\left|\delbar_J u_i\right|_{g^{cl}} \leq \frac{2\sqrt{R_\LL}\cdot \ep_{(F,J)}}{L_i} \to 0
\ee
uniformly on compact subsets of $z$'s with $|z - z_i| \leq \delta_i L_i$
as $i \to \infty$.
Combining \eqref{eq:eizi}, \eqref{eq:delbarto0}, we derive a subsequence of
$\widetilde u_i$ so that it converges to a $J_\infty$-holomorphic map
$$
v: \C \subset \C \cup \{\infty\} \cong S^2  \to M
$$
or
$$
v: (\H, \del \H) \cong (D^2 \setminus \{1\}, \del D^2 \setminus \{1\}) \to (M, L) \quad \mbox{or } \, (M,L_0)
$$
such that $|dv(0)| \geq 1$ and $\operatorname{Im} v \cap M^{\leq R} \neq \emptyset$.

Therefore one of the following alternatives should hold:
\begin{enumerate}
\item $\operatorname{Im} v$ is bounded and so $v$ smoothly extends either to $S^2$
or $(D^2, \del D^2)$
\item $\operatorname{Im} v$ is unbounded
\end{enumerate}

The first case is ruled out by the exactness hypothesis on $(L,L_0;M)$.
For the second case, due to $\operatorname{Im} v \cap M^{\leq R} \neq \emptyset$
and the are bound
\bea\label{eq:Ev}
\operatorname{Area}_{g^{cl}}(v) & \leq & \operatorname{Area}_{g_{(\psi,J)}}(v) \leq
\lim_{i \to \infty} E_{g_{(\psi,J^\chi)}}(\widetilde u_i) + 2R_\LL\, \ep^2_{(F,J)}\nonumber \\
& \leq & 2C_1 + 2 R_\LL\, \ep^2_{(F,J)},
\eea
we can apply the Courant-Lebesgue Lemma to find a sequence of discs $D^2(R_i) \subset \C$
or $D^2(R_i)\subset \H$ with radi $R_i \to \infty$ such that
\be\label{eq:lengthto0}
\operatorname{length}(v|_{\del D^2(R_i)}) \to 0.
\ee
(See \cite{oh:removal} for the case with Lagrangian boundary).
Then this combined with the unboundedness of $v$ and the monotonicity formula
implies that the $g^{cl}$-area of $v$ must be infinite,
which clearly contradicts to the area bound \eqref{eq:Ev}.

\begin{rem} In fact, in the current exact context, we can derive a contradiction
immediately from \eqref{eq:Ev}-\eqref{eq:lengthto0} as follows: By the finiteness of
the area, we have the identity
\beastar
\operatorname{Area}_{g^{cl}}(v) & = & \lim_{i \to \infty} \int_{D^2(R_i)} v^*\omega \\
& = & \lim_{i\to \infty} \int_{\del D^2(R_i)} v^*\alpha \leq
|\alpha|_{C^0} \, \operatorname{length}(v|_{\del D^2(R_i)}) \to 0.
\eeastar
Here $|\alpha|_{C^0}$ denotes the $C^0$-norm of $\alpha$. It is bounded
since $(M, \widehat \omega = d\alpha)$ is geometrically bounded.
This gives rise to a contradiction since $v$ is non-constant.

We prefer to give the above more general proof since the argument using the
monotonicity together with Courant-Lebesgue lemma, which is based on the
unboundedness of the image of $v$, can be applied to non-exact context.
\end{rem}

Next we consider the case $R_i:= r(u_i(z_i)) \to \infty$.
We split our consideration into the two cases, one where $R_i - (2\delta_i L_i+1) \to +\infty$ and the
other where $R_i - (2\delta_i L_i+1) \leq C_3$ for some $C_3$.

If $R_i - (2 \delta_i L_i+1) \to \infty$, we have
$$
\widetilde u(D^2(\delta_i L_i)) = u_i(D_{z_i}(\delta_i)) \subset M^{R_\LL+1} \cong [\ln(R_\LL+1), \infty) \times Q \subset \R \times Q.
$$
Therefore we may regard $\widetilde u_i$ as a map from $D^2(\delta_i L_i)$ to $\R \times Q$.
We also have
\be\label{eq:derivbound}
|d\widetilde u_i(z)|_{g^{cl}} \leq 2
\ee
for all $z \in D^2(\delta_i L_i)$. Write $\widetilde u_i = (a_i, w_i)$ where
$a_i = s\circ \widetilde u_i$ and $w_i = \pi \circ \widetilde u_i$
and consider the translated map
$$
\overline u_i(z) = (a_i(z) - \ln R_i, w_i).
$$
Since this translation does not change the gradient and energy bound
(with respect to $g^{cl}$), we still have the same bounds for $\overline u_i$ as
$\widetilde u_i$.

Again we derive a subsequence of
$\widetilde u_i$ so that it converges to a $J_\infty$-holomorphic map
$$
v: \C \subset \C \cup \{\infty\} \cong S^2  \to \R \times Q
$$
or
$$
v: (\H, \del \H) \cong (D^2 \setminus \{1\}, \del D^2 \setminus \{1\}) \to (\R \times Q, \R \times \del_\infty L)
\quad \mbox{or } \, (\R \times Q, \R \times \del_\infty L_0)
$$
such that $|dv(0)| \geq 1$ and $\operatorname{Im} v \cap M^{\leq R} \neq \emptyset$,
and $E_{g^{cl}}(v) \leq 2C_1 +1$.
We obtain a contradiction as before.
\par
Next we consider the case where $R_i - 2\delta_i L_i+1 \leq C_3$.
We set
$$
R_1 = \max\{R_\LL, \sup_{i} R_i - (2\delta_i L_i+1)\}.
$$
\begin{lem} Let $M_i = \frac{R_i - R_1}{2}$. Then we have
$$
\widetilde u_i(D^2(M_i)) \subset M^{\geq R_\LL} \cong [\ln R_\LL, \infty) \times Q
$$
where the last product is in terms of cylindrical coordinates $(s,y)$.
\end{lem}
\begin{proof} As before we estimate $s(\widetilde u_i(z)) - \ln R_i$ with $R_i = r(\widetilde u_i(z_i))$
on $z \in D^2(M_i)$
$$
|s(\widetilde u_i(z)) - \ln R_i| \leq 2M_i = \ln R_i - \ln R_1
$$
and hence $s(\widetilde u_i(z)) \geq \ln R_1$ for all $z \in D^2(M_i)$. Since $R_1 \geq R_\LL$ by definition,
we have finished the proof.
\end{proof}

We remark that $M_i = \frac{R_i - R_1}{2} \to \infty$.
Therefore we can now repeat the above process used for the case
$R_i - (2\delta_i L_i+1) \to \infty$ to the map $\widetilde u_i$ restricted on $D^2(M_i)$
and get a contradiction.
\par
Therefore \eqref{eq:maxderivbound} should hold. This finishes the proof.

\end{proof}

\section{$C^0$ bound}
\label{sec:C0bound}

In this section we will use the energy bound and the derivative bound to prove the a priori
$C^0$ bound for a \emph{$C^1$-small isotopy}. We will quantify the $C^1$-smallness
of $F$ that we need in terms of
$
\ep_{(F,J)}.
$
Let $a_F^{\infty} = |f|_\infty + |df|_\infty $ for $F(r,y) = r f(y)$ at
the end as before. In particular $\ep_{(F,J)}^{out} \leq a_F^\infty$ is bounded
for any given such $F$. We also define
$$
R_{\operatorname{Int} \LL}:= \inf\{R \mid \operatorname{Int} \LL \subset M^{\leq R}\}.
$$
The following is the main $C^0$ bound we prove in this section.

\begin{thm}\label{thm:C0boundforsmall} Let $(M,\omega)$ be Weinstein and $L, \, L_0$ be
conic with their thresholds $R_L, \, R_{L_0}$ respectively.
Let $\CL = \{L_s\}$ be a Hamiltonian isotopy generated by a one-parameter family of
conic Hamiltonians $F = \{F_s\}$ such that
\begin{enumerate}
\item $(L_0,L)$ and $(L_1,L)$ are transversal pairs.
\item $\operatorname{Int} \LL=\cup_{s \in [0,1]}L \cap L_s$ is compact.
\end{enumerate}
Then there exists a sufficiently small $\ep_0 = \ep_0(R_\LL))> 0$ such that for
any conic Hamiltonian $F$ with $\ep^{out}_F \leq \ep_0$ we have
\be\label{eq:maximum-ru}
\operatorname{Im} u \subset M^{\leq R_\LL}.
\ee
We set $C_6 = R_\LL$.
\end{thm}
\begin{proof}
We prove this by contradiction. Suppose to the contrary that there exists a sequence of
conic Hamiltonians $F_i$ such that
\be\label{eq:eFiout}
\ep_{(F_i,J)}^{out} \to 0
\ee
and a map $u_i: \R \times [0,1] \to M$ satisfying
\be\label{eq:dudtauXFdudt}
\dudtau - X_{F_i}^\chi(u) + J^\chi\left(\dudt-X_{F_t}(u)\right) = 0
\ee
but $\max (r\circ u_i) > R_\LL$.
By the energy bound we obtained in \eqref{eq:energy}-\eqref{eq:boundC1}, we have the energy
\be\label{eq:energybound}
E_{(g_{(\psi,J^\chi)},F_i)}(u_i) \leq 2w_\CA(L,L_0) +1 =: C_4
\ee
when $i$ is sufficiently large so that $6 R_\LL \ep_{(F_i,J)}^2 \leq 1$.
We emphasize that $C_4$ is independent of $i$'s.

Since $F_i^\chi \equiv 0$ for $\tau \geq 1$ or for $\tau \leq 0$,
$u_i$ satisfies the genuine Cauchy-Riemann equation
$$
\dudtau + J^\chi \dudt = 0
$$
on $[1,\infty) \times [0,1] \cup (-\infty,0] \times [0,1]$. Therefore
by the maximum principle and strong maximum principle, the maximum must
be achieved at some $z_i = (\tau_i,t_i)$ with $0 \leq \tau_i \leq 1$
whenever $\operatorname{Im} u \not \subset M^{\leq R_\LL}$.
We also note that $[0,1] \times [0,1]$ is compact and $F_i^\chi \to 0$.

We first prove that $r(u(z_i))$ is bounded. If not,
by choosing a subsequence,
we assume that $(\tau_i,t_i) \to (\tau_\infty,t_\infty)$
as $i \to \infty$. Then we have $(\tau_\infty,t_\infty) \in [0,1] \times [0,1]$.
We note that $[0,1] \times [0,1]$ is compact and so \eqref{eq:dudtauXFdudt}
converges to the Cauchy-Riemann equation.

We consider the restriction
$$
u_i: [-1,2] \times [0,1] \to M.
$$
By the derivative bound $|du_i| \leq C_2$ and the assumption $R_i:= r(u_i(z_i)) \to \infty$,
it follows that we have
$$
u_i([-1,2] \times [0,1]) \subset M^{R_i-C_5 \leq r\leq R_i + C_5}
\subset [R_i -C_5, \infty) \times Q \subset \R \times Q
$$
for some constant $C_5$ independent of $i$'s.
Then since $u_i$ is $J$-holomorphic outside $[0,1] \times [0,1]$ and have the derivative bound,
$$
u_i: u_i^{-1}\left(M^{\leq R_i-C_5}\right) \to M^{\leq R_i-C_5}
$$
defines a proper $J_0$-holomorphic curve.
We note that $\Sigma_i:= u_i^{-1}\left(M^{\leq R_i -C_5}\right)$
has boundary components consisting of 3 types
\beastar
\del \Sigma_i & = & (\del \Sigma_i \cap \R \times \{0\}) \bigcup (\del \Sigma_i \cap \R \times \{1\}) \\
& {} & \quad \bigcup \del (\Sigma_i \cap u_i^{-1}(M^{R_i-C_5})).
\eeastar
Using the bounded geometry of $(M,L,L_0)$ and $u(\pm \infty) \in M^{\leq \max\{R_L,R_{L_0}\}}$,
we can apply the monotonicity formula (see Proposition 4.3.1 \cite{sikorav} for the precise
formulation of the required monotonicity formula) to conclude
$\operatorname{Area}(u_i(u_i^{-1}(M^{\leq R_i - C_5})))$ must be infinite, which
contradicts to the above area bound. This finishes the proof of
boundedness of the image.

Knowing the boundedness of the image, we choose a sufficiently large $R_3 > 0$ such that
$$
\operatorname{Im } u \subset M^{\leq R_3}
$$
for all finite energy solution $u$ of \eqref{eq:dudtauXFdudt} with bounded image.

Let $\delta > 0$ and suppose that $r(u_i(z_i)) > R_\LL+\delta$
for all $i$ with $z_i(\tau_i,t_i) \in [0,1]\times [0,1]$
Since the images of $u$'s all lie in a common compact subset $M^{\leq R_3}$,
and since $[0,1] \times [0,1] \subset \R \times [0,1]$ is a compact subset,
we can use the energy bound and apply Gromov-Floer compactness theorem to the sequence
$u_i$ to extract a $J_0$-holomorphic curve $u$ satisfying
$$
u(\pm\infty) \in L\cap L_0, \quad r(u(\tau_0,t_0)) \geq R_{\LL} +\delta\quad (\tau_0,t_0) \in [0,1] \times [0,1]
$$
as $i \to \infty$.
Here the second condition comes from the hypothesis $r(u(z_i)) \geq R_\LL + \delta$. However on
$M^{\geq R_\LL}$, both $L$ and $L_0$ are conic and hence we can apply the maximum
principle and strong maximum principle to derive that the latter is impossible
since $u(\pm\infty) \in M^{\leq R_{\operatorname{Int}\LL}} \subset M^{\leq R_\LL}$. This proves
$$
\operatorname{Im } u \subset M^{\leq R_\LL + \delta}.
$$
Since $\delta$ is arbitrary, we have finished the proof.
\end{proof}

\begin{rem} If we allow the constant $C_6$ to depend on the second order behavior of
$L, \, L_0$, i.e., on the norms of the second fundamental forms of $L, \, L_0$
and combine the monotonicity formula, Proposition \ref{prop:monotonicity},
we can slightly improve the above $C^0$-estimate as follows.

We fix the constant $r_0> 0$ appearing in Proposition \ref{prop:monotonicity}.
Then we can find at least $\frac{R_3-R_{\operatorname{Int}\LL}}{2r_0}$
disjoint balls $B(x,r_0)$ with $x \in \operatorname{Im}u$. Therefore the above
monotonicity formula implies
$$
C_5 r_0^2 \cdot \frac{R_3-R_{\operatorname{Int}\LL}}{2r_0} \leq C_4
$$
i.e.,
we obtain
\be\label{eq:R3estimate}
R_3 \leq R_{\operatorname{Int}\LL} + \frac{2C_4}{C_5 r_0}.
\ee
Therefore we can take $C_6$ to be
$$
\min\left\{R_{\operatorname{Int}\LL} + \frac{2C_4}{C_5 r_0}, R_\LL\right\}.
$$
\end{rem}

\section{Wrap-up of construction of the chain map}
\label{sec:wrapup}

Now we go back to the equation \eqref{eq:CRchain}. Take a partition
$$
P: 0 = s_0 < s_1 < \cdots < s_K = 1
$$
and denote $\delta = \operatorname{mesh}P$. It follows that if $\delta$ is sufficiently small
we can connect $L_{s_k}$ to $L_{s_{k+1}}$ by the Hamiltonian isotopy
generated by a conic Hamiltonian $F_k$ such
$$
\ep_{(F_k, \phi_{F_k}^*J)}
$$
is sufficiently small so that we have the a priori $C^0$ bound: The equation \eqref{eq:CRchain}
is transferred \eqref{eq:XFXHchi} with $F = F_k$ and $J = \phi_{F_k}^*J$.
By letting the partition finer, we can apply the energy and $C^0$ bound for the
family $(F,J) = (F_k, \phi_{F_k}^*J)$ and hence give rise to the chain
map $h_{F_k}: CF(L_{s_k},L) \to CF(L_{s_{k+1}},L)$, \emph{whenever} the
partition $P$ is chosen so that $(L_{s_k},L)$ and $(L_{s_{k+1}},L)$ are
Floer regular in the sense that all the relevant moduli spaces are Fredholm
regular.
\par
In fact, Theorem 4.6 \cite{oh:jkms} states that for the Hamiltonian isotopy generated by
a generic $F$ so that there exists a partition $P$ so that the sub-homotopies
from $s_k$ to $s_{k+1}$ for all $k=0, \cdots, K-1$ are all Floer regular
so that the  chain map $h_{F_k}: CF(L_{s_k},L) \to CF(L_{s_{k+1}},L)$
can be constructed by considering the Floer chain map equation associated
to the sub-homotopies $\{F_s\}_{s_k \leq s\leq s_{k+1}}$.
This allows us to use the sub-homotopy of
$\phi_F^s$ for $s_k \leq s \leq s_{k+1}$ for the given Hamiltonian isotopy.
While the construction
explicitly depends on the choice of $J$, the argument used in \cite{oh:jdg}
proves that this map, \emph{including filtrations}, is independent of the
choice of $J$, and so we drop $J$-dependence from this homomorphism.

We then take the adiabatic chain map
$$
h_{\CL;P}:= h_{F_K} \circ h_{F_{K-1}} \circ\cdots \circ h_{F_0} :
CF(L_0,L) \to CF(L_1,L)
$$
and take its homology.
By a simple adiabatic argument employed in \cite{milinko-oh}, one can prove that
this chain map does not depend on the choice of the partition $P$ as long as $\operatorname{mesh}P$ is
sufficiently small, and hence depends only on the given isotopy $\CL = \{L_s\}$.
This shows that the chain map is canonical and depends only on the isotopy $\CL$
(or more precisely on the Hamiltonian $F$ generating it).

\section{The case of asymptotically conic Lagrangian submanifolds}
\label{sec:asympconic}

In this section, we explain how we can modify our construction to accommodate
$C^2$-asymptotically conic Lagrangian submanifolds.

\subsection{Approximation by conic Lagrangian submanifolds}
\label{subsec:normalform}

We first give a description of an asymptotically conic Lagrangian submanifold
in terms of its asymptototic boundary $N = \del_\infty L$.
Let $R_L$ be any constant given in Definition \ref{defn:asympconic}
and $N = \del_\infty L$. By Darboux theorem, there is neighborhood
$U$ of $N$ in $Q$ and $V$ of the zero section $o_N$ of the one-jet bundle
$J^1(N)$ such that we have a contactomorphism
$$
\Psi : U \subset Q \to V \subset J^1(N)
$$
such that $\lambda = \Psi^*(\theta - dz)$ where $\theta = pdq$ is the
Liouville one form on $T^*N$ (see Appendix \cite{arnold:book}) and
restricts to the identity on $N$ via the identification $N \cong o_N$.
Under this isomorphism, the Lagrangian submanifold $[R_L,\infty) \times N$
corresponds to $[R_L,\infty) \times \{0\} \times o_N \subset
[R_L,\infty) \times \R \times T^*N$.

By the $C^1$ uniform convergence of $L(R)$ to $N = \del_\infty L $,
we can write $L \cap [R_L,\infty) \times Q$
\be\label{eq:rzqp}
\{(r,z,q,p)  \mid p = \alpha(r,q), \, z = h(r,q), \, q \in N, \, r \in [R_L,\infty)\}.
\ee
\begin{prop}\label{prop:Lag} The subset \eqref{eq:rzqp} is Lagrangian
in $\R_+ \times \R \times T^*N$ with respect to the symplectic form
$$
d(r(pdq - dz)) = dr \wedge (pdq -dz) - r(dq \wedge dp)
$$
if and only if $\alpha$ and $h$ satisfies
\be\label{eq:alphah}
\alpha(r,q) = \frac{1}{r}\left( \int_{R_L}^r d_N h(a,q)\, da + \beta(q)\right)
\ee
where $\beta$ is a closed one form on $N$.
\end{prop}
\begin{proof} The symplectic form $\omega = d(r \pi^*\lambda)$ is given
by
$$
d(r\pi^*(pdq - dz)) = dr \wedge \pi^*(pdq -dz) + r d(\pi^*(pdq -dz)).
$$
Therefore we have
\beastar
i_L^*(d(r\pi^*(pdq - dz))) & = & dr \wedge \alpha -d_Nh + r dr \wedge \frac{\del \alpha}{\del  r}
- r d_N\alpha \wedge dq \\
& = & dr\wedge \left(\alpha + r \frac{\del \alpha}{\del  r} -d_Nh\right) - r d_N \alpha \wedge dq.
\eeastar
Hence $i_L^*(d(r\pi^*(pdq - dz))) = 0$ if and only if
$$
\alpha + r \frac{\del \alpha}{\del  r} -d_Nh = 0, \quad d_N \alpha \wedge dq = 0.
$$
The second equation implies that $d_N \alpha = 0$, i.e., $\alpha$ is a $r$-family of
closed one-forms on $N$. On the other hand we can write the first equation as
$$
\frac{\del(r\alpha)}{\del r} = d_N h.
$$
By integrating the equation over $r$, we obtain
$$
r\alpha(r,q) = \int_{R_L}^r d_N h(a,q)\, da + \beta(q)
$$
where $\beta$ is a closed one-form on $N$. Therefore we obtain
$$
\alpha(r,q) = \frac{1}{r} \left(\int_{R_L}^r d_N h(a,q)\, da + \beta(q)\right).
$$
This finishes the proof.
\end{proof}

\subsection{Canonical adiabatic homomorphism}
\label{subsec:canonicalhomo}
Assume that $L$ is $C^2$-asymptotically conic. We note
$$
m_\gamma^{-1}(L \cap M^{\geq \gamma R_L})
= \{(r,z,q,p) \mid p = m_\gamma^*\alpha, \, z = m_\gamma^*h, \, q \in N, \, r \in [R_L,\infty) \}.
$$
Therefore the $C^2$-asymptotic conic condition becomes
$$
\lim_{\gamma \to \infty} m_\gamma^*\alpha = 0 = \lim_{\gamma \to \infty} m_\gamma^* h
$$
on $[R_L,\infty) \times N$ uniformly in $C^2$-topology.

Now we fix a monotonically decreasing cut-off function $\rho: [1,\infty) \to [0,1]$
such that $\rho(r) \equiv 1$ for $r \in [1,R_L]$ and $\rho(r)\equiv 0$ for $r \geq R_L+1$.
Then we define a two-parameter family of functions $\rho_{R,K}: [1, \infty) \to [0,1]$ given
by
$$
\rho_{R,K}(r) = \rho\left(\frac{r - R}{K}\right)
$$
for all $R \geq 0$ and $K \geq 1$. Then we define a two-parameter family of
Lagrangian submanifolds defined by the pair $(\alpha_{R,K},h_{R,K})$ where
\beastar
h_{R,K}(r,q) & = & \rho_{R,K}(r) h(r,q), \\
\alpha_{R,K} & = &
\frac{1}{r} \left(\int_{R_L}^r d_N \rho_{R,K}(r) h(a,\cdot)\, da + \beta\right).
\eeastar
The following is easy to see whose proof we leave to the readers.

\begin{prop} Let $R \geq 0$. Then $L_{R,K}$ is a conic Lagrangian
submanifold with its threshold given by $R_L + R +K$ and converges to
$L$ in $C^\infty$-topology as $R \to \infty$.
\end{prop}

Note that under this approximation, the intersection set
$$
\operatorname{Int}\LL= \bigcup_{s \in [0,1]} L \cap L_s
$$
is unchanged and all the norms of second fundamental forms of
the family $L_{R,K}$ are uniformly bounded over $R, \, K$ by the $C^2$-
asymptotically conic property of $L, \, L_0$'s. From now on, we will fix $K=1$ and
vary $R$ only. For any given pair $R' > R$, we consider the isotopy
$$
\LL^\rho =\{L_s\}_{s \in [0,1]}, \quad L_s = L_{(1-s)R + s R',1}
$$
where $\rho$ is the cut-off function used in the definition of the family
$L_{R,1}$. Then the following lemma can be proved by the argument from
\cite{oh:mathz,oh:jkms}.

\begin{lem} For a generic choice of $\rho$, there exists a dense
subset $I_\rho \subset [0,1]$ such that $0, \, 1 \in I_\rho$
and for any given pair $(s,s') \in I_\rho$
the subhomotopy  $\LL|_{[s,s']}$ is Floer-regular in the sense that
the Floer chain map $h_{ss'}=h_{LL|_{[s,s']}}: CF(L_s,L) \to CF(L_{s'},L)$
can be defined.
\end{lem}

We fix such a $\rho$ and define the isotopy $\LL^\rho$. Then we
choose a partition
$$
P^\rho: 0=s_0 < s_1 < \cdots < s_N =1
$$
so that $\operatorname{mesh}P$ is so small that the corresponding
$\ep_0$ satisfies Theorem \ref{thm:C0boundforsmall}, and construct the corresponding
chain map $h_{ss'}: CF(L_s,L) \to CF(L_{s'},L)$.

The following proposition proves that the Floer cohomology $HF(L_0,L)$
can be described as the inverse limit of $HF(L_{0,R},L_R)$
as $R \to \infty$.
This easy to prove from our discussions whose proof we omit noting that
the family $\{(L_{0,R},L_R)\}$ can be made into a \emph{monotone} family
so that the corresponding chain map can be easily constructed by applying
the maximum principle. We omit its proof and
refer readers to \cite{kasturi-oh1} for the details of the proof employed in a similar
context of approximations of the conormals of open sets.

\begin{prop} Denote by $\iota_{RR'}: CF(L_{0,R'}, L_R') \to CF(L_{0,R},L_{R})$ be the
canonical map obtained from the identification $L_R \cap L_{0,R} = L_{R'} \cap L_{0,R'}$. Then
\begin{enumerate}
\item $\iota_{RR'}$ is a chain map
for any $R' > R$ which respects the diagram
\be\label{eq:squareRR'}
\xymatrix{
CF(L_{0,R'}, L_{R'}) \ar[r]_{\iota_{RR'}} \ar[d]_{\del_{R'}} & CF(L_{0,R},L_{R}) \ar[d]_{\del_{R}} \\
CF(L_{0,R'}, L_{R'}) \ar[r]_{\iota_{RR'}} &  CF(L_{0,R},L_{R})
}
\ee
and hence $\{\iota_{RR'}: CF(L_{0,R'}, L_{R'}) \to CF(L_{0,R},L_{R})\}$ defines an
inverse system of chain maps.
\item There exists a canonical homomorphism
$$
(\iota_R)_*: CF(L_0,L) \to CF(L_{0,R},L_{R})
$$
such that the following diagram
$$
\xymatrix{ HF(L_{0,R'}, L_{R'}) \ar[rr]^{(\iota_{RR'})_*} &{} & HF(L_{0,R},L_{R})\\
{} & HF(L_0,L)\ar[ul]^{(\iota_{R'})_*}\ar[ur]_{(\iota_R)_*} &{} }
$$
commutes.
\item The above chain map $(\iota_R)_*$ induces an isomorphism in homology.
\end{enumerate}
\end{prop}

Once we have established this description of $HF(L_0,L)$ for
any asymptotically conic pair $(L_0,L)$, we can define the natural
isomorphism
$$
h_{\LL}: HF(L_0,L) \to HF(L_1,L); \quad \LL = \{L_s\}_{s \in [0,1]}
$$
as the inverse limit of the homomorphisms
$$
h_{\LL_R}: HF(L_{0,R},L_R) \to HF(L_{1,R},L_R)
$$
as $R \to \infty$: One can easily prove the natural adiabatic chain map
$$
h_{\LL_R}: CF(L_{0,R},L_{R}) \to CF(L_{1,R},L_R)
$$
is a chain isomorphism. Therefore the inverse limit
$$
\lim_{\longleftarrow} h_{\LL_R}
$$
is an isomorphism. Then composition map
$$
(\lim_{\longleftarrow}(\iota_{1,R})_*)\circ(\lim_{\longleftarrow}
h_{\LL_R})\circ (\lim_{\longleftarrow}(\iota_{0,R})_*)^{-1} : HF(L_0,L) \to HF(L_1,L)
$$
becomes an isomorphism. We denote this isomorphism by $h_{\LL}$.
We would like to point out that we do not have the relevant chain map
$CF(L_0,L) \to CF(L_1,L)$ in which sense the current situation is
similar to that of \cite{kasturi-oh1,kasturi-oh2}.

\section{Discussions of the non-exact case and the $A_\infty$ chain map case}
\label{sec:nonexact}

For the non-exact $(M,\omega)$ with cylindrical ends, usual consideration of bubbling
spheres or discs should be taken into consideration and hence
requires the whole consideration of $A_\infty$ structures to do
the Lagrangian Floer theory as performed in \cite{fooo:book}.

The only thing we have to make sure beside this machinery before carrying
out the Floer theory is to obtain the a priori energy and $C^0$ bounds.
It turns out that we can prove the following  energy bound
by combining the idea of the proof in this paper with the ones in
\cite{oh:mrl}, \cite{abou-seidel}.

\begin{thm} Let $(M,\omega)$ be a symplectic manifold convex at infinity.
Assume that $J$ is of contact type and $F$ is conic at end.
Consider a pair of anchored Lagrangian
submanifolds $(L, L_0)$ in the sense of \cite{fooo:anchored}
such that $L, \, L_0$ are either compact or have
conic ends or their mixture. Fix a homotopy class $A \in \pi_2(p^-,p^+;L_0,L;M)$
for $p^\pm \in L\cap L_0$.
Let $\CL = \{L_s\}$ be a Hamiltonian isotopy generated by $F = \{F_s\}$
such that
\begin{enumerate}
\item $(L_0,L)$ and $(L_1,L)$ are transversal pairs.
\item $\cup_{s \in [0,1]}L \cap L_s$ is compact.
\end{enumerate}
Then there exists a constant $C_6$ for
$i=1, \, 2$ such that every finite energy solution $u$ of
\eqref{eq:FHCR} with
$$
u(-\infty) = p^-, \, u(\infty)=p^+, \quad [u] = A,
\,
$$
satisfies
\be\label{eq:energyC6} E_{g_{(\psi,J)}}(u) \leq C_6
\ee
where $C_6$ depends only on $A$ and $(F,J)$ but independent of $u$.
\end{thm}
In fact the scheme of our proof of the
energy bound \eqref{eq:energy} does not require the isotopy $\CL$
to be Hamiltonian but only to be a symplectic isotopy. In other words,
the Hamiltonian vector field $X_F$ in the perturbed Cauchy-Riemann equation
\eqref{eq:XFXHchi} by any symplectic (or locally Hamiltonian) vector fields.

Once we have established this energy estimate, essentially the same proof of
$C^0$-bound applies to the general symplectic manifolds convex at infinity.

Taking this process for granted, which is reasonable to accept,
the rest of the Floer theory will be a straightforward
repetition of those carried out in
\cite{fooo:book} for the current case of chain maps in homology, and
in \cite{seidel:book} for the setting of $A_\infty$ maps.
The outcome will be the analog to the theory developed in
\cite{fooo:book} for unwrapped Fukaya category generated by
the compact and (asymptotically) conic Lagrangian submanifolds, not necessarily exact,
on symplectic manifolds Weinstein at infinity.

We will come back to this generalization elsewhere.

\end{document}